\input ./style/arxiv-general.cfg
\documentclass[aap,MSNbibl,seceqn,dvips]{arximspdf}
\makeatletter
   \@ifpackageloaded{graphicx}{}{\usepackage{graphicx}}
\makeatother

\doi{10.1214/14-AAP1091}
\volume{26}
\issue{1}
\pubyear{2016}
\firstpage{305}
\lastpage{327}
\docsubty{FLA}

\makeatletter
\newcommand{\ds}{\displaystyle}
\newcommand{\rrvert}{\vert}
\newcommand{\llvert}{\vert}
\newtheorem{lemma}{Lemma}[section]
\newtheorem{theorem}{Theorem}
\newtheorem{proposition}{Proposition}[section]
\makeatother

\begin{document}
\begin{frontmatter}

\title{The densest subgraph problem in sparse random~graphs\thanksref{T1}}
\runtitle{Densest subgraphs in sparse graphs}

\begin{aug}
\author[A]{\fnms{Venkat}~\snm{Anantharam}\ead[label=e1]{ananth@berkeley.edu}}
\and
\author[B]{\fnms{Justin}~\snm{Salez}\corref{}\ead[label=e2]{justin.salez@univ-paris-diderot.fr}}
\runauthor{V. Anantharam and J. Salez}
\affiliation{University of California, Berkeley and Universit\'e Paris Diderot (LPMA)}
\address[A]{University of California, Berkeley\\
Department of Electrical Engineering\\
\quad and Computer Sciences\\
253 Cory Hall\\
Berkeley, California 94720-1770\\
USA\\
\printead{e1}}
\address[B]{Universit\'e Paris Diderot\\
UFR de Math\'ematiques\\
5 rue Thomas Mann\\
75205 Paris CEDEX 13\\
France\\
\printead{e2}}
\end{aug}
\thankstext{T1}{Supported by ARO MURI Grant W911NF-08-1-0233,
\textit{Tools for the Analysis and Design of Complex Multi-Scale
Networks}, NSF Grant CNS-0910702, and
NSF Science \& Technology Center Grant CCF-0939370, \textit{Science of
Information}.}

%
\received{\smonth{1} \syear{2014}}
%
\revised{\smonth{7} \syear{2014}}

%
\begin{abstract}
We determine the asymptotic behavior of the maximum subgraph density of
large random graphs with a prescribed degree sequence. The result
applies in particular to the Erd\H{o}s--R\'enyi model, where it settles
a conjecture
of Hajek [\textit{IEEE Trans. Inform. Theory} \textbf{36} (1990) 1398--1414]. Our proof consists in extending the notion of balanced
loads from finite graphs to their local weak limits, using
unimodularity. This is a new illustration of the objective method
described by Aldous and Steele [In \textit{Probability on Discrete Structures}  (2004)
1--72 Springer].
\end{abstract}

%
\end{frontmatter}

\section{Introduction}

Let $G=(V,E)$ be a finite, simple, undirected graph. Write $\vec{E}$
for the set of oriented edges, formed by replacing each edge $\{i,j\}
\in E$ with the two oriented edges $(i,j)$ and $(j,i)$. An \textit{allocation} on $G$ is a map $\theta\dvtx\vec{E}\to[0,1]$ satisfying
$\theta(i,j)+\theta(j,i)=1$ for every $\{i,j\}\in E$. The \textit{load}
induced by $\theta$ at a vertex $o\in V$ is
\[
\partial\theta(o)  :=  \sum_{i\sim o}
\theta(i,o),
\]
where $\sim$ denotes adjacency in $G$. $\theta$ is \textit{balanced} if
for every $(i,j)\in\vec{E}$,
\begin{equation}
\label{eqbalance}
\partial\theta(i)<\partial\theta(j) \quad \Longrightarrow\quad
\theta(i,j)=0.
\end{equation}
Intuitively, one may think of each edge as carrying a unit amount of
load, which has to be distributed over its end-points in such a way
that the total load is as balanced as possible across the graph. In
that respect, (\ref{eqbalance}) is a \textit{local optimality}
criterion: modifying the allocation along a single edge cannot further
reduce the load imbalance between its end-points. This condition
happens to guarantee \textit{global optimality} in a very strong sense.
Specifically, the following conditions are equivalent (see
\cite{Haj90}):
\begin{longlist}[(iii)]
\item[(i)] $\theta$ is balanced.
\item[(ii)] $\theta$ minimizes $\sum_{o\in V}f (\partial\theta
(o)
)$, for some strictly convex $f\dvtx[0,\infty)\to\mathbb{R}$.
\item[(iii)] $\theta$ minimizes $\sum_{o\in V}f (\partial\theta
(o))$, for every convex $f\dvtx[0,\infty)\to\mathbb{R}$.
\end{longlist}

In particular, balanced allocations exist on $G$ and they all induce
the same loads $\partial\theta\dvtx V\to[0,\infty)$. The balanced load
$\partial\theta(o)$ induced at a vertex $o\in V$ has a remarkable
graph-theoretical interpretation: it measures the \textit{local density}
of $G$ at $o$. Specifically, it was shown in \cite{Haj90} that the
vertices receiving the highest load solve the {classical densest
subgraph problem} on $G$: the value $\max\partial\theta$ coincides with
the \textit{maximum subgraph density} of $G$,
\[
\varrho(G)  :=  \max_{\varnothing\subsetneq H\subseteq V}\frac{|E[H]|}{|H|},
\]
and the set $H=\operatorname{argmax}\partial\theta$ is precisely the
largest set achieving this maximum. Here, $E[H]\subseteq E$ naturally
denotes the set of edges with both end-points in $H$. This surprising
connection with a well-known and important graph parameter justifies a
deeper study of balanced loads in large graphs. It is convenient to
encode the loads induced by a balanced allocation on $G$ into a
probability measure on $\mathbb{R}$, called the \textit{empirical load
distribution} of $G$:
\[
\mathcal{L}_G  =  \frac{1}{|V|}\sum
_{o\in V} \delta_{\partial
\theta(o)}.
\]

Motivated by the above connection, Hajek \cite{Haj90} studied the
asymptotic behavior of $\mathcal{L}_{G}$ on the classical \textit{Erd\H
{o}s--R\'enyi model},
where the graph $G=G_n$ is chosen uniformly at random among all graphs
with $m=\lfloor\alpha n\rfloor$ edges on $V=\{1,\ldots,n\}$. In the
regime where the density parameter $\alpha\geq0$ is kept fixed while
$n\to\infty$, he conjectured that $\mathcal{L}_{G_n}$ should
concentrate around
a deterministic probability measure $\mathcal{L}\in\mathcal
{P}(\mathbb{R})$ and that
\[
\varrho(G_n)\mathop{\longrightarrow}_{n\to\infty}^{\mathbb{P}}
\varrho :=\sup\bigl\{t\in\mathbb{R}\dvtx\mathcal{L} \bigl([t,+\infty) \bigr)>0
\bigr\}.
\]
Using a nonrigorous analogy with the case of finite trees, Hajek even
proposed a description of $\mathcal{L}$ and $\varrho$ in terms of the
solutions
to a distributional fixed-point equation which will be given later.
In this paper, we establish this conjecture together with its analogue
for various other sparse random graphs, using the unifying framework of
\textit{local weak convergence}.

\section{Local weak convergence}
\label{seclwc}
This section gives a brief account of the theory of local weak
convergence. For more details, we refer to the seminal paper \cite
{BenScr01} and to the surveys \cite{AldSte04,AldLyo07}.

\subsection*{Rooted graphs}
A \textit{rooted graph} $(G,o)$ is a graph
$G=(V,E)$ together with a distinguished vertex $o\in V$, called the
\textit{root}. We let $\mathcal{G}_{\star}$ denote the set of all
locally finite
connected rooted graphs considered up to \textit{rooted isomorphism},
that is, $(G,o)\equiv(G',o')$ if there exists a bijection $\gamma\dvtx
V\to V'$ that preserves roots ($\gamma(o)=o'$) and
adjacency ($\{i,j\}\in E\Longleftrightarrow\{\gamma(i),\gamma(j)\}
\in
E'$). We write $[G,o]_h$ for the (finite) rooted subgraph induced by
the vertices lying at graph-distance at
most $h\in\mathbb{N}$ from $o$. The distance
\[
\operatorname{DIST} \bigl((G,o),\bigl(G',o'\bigr) \bigr):=
\frac{1}{1+r} \qquad \mbox{where } r=\sup \bigl\{h\in\mathbb{N}
\dvtx [G,o]_h\equiv\bigl[G',o'
\bigr]_h \bigr\},
\]
turns $\mathcal{G}_{\star}$ into a complete separable metric space;
see \cite{AldLyo07}.

\subsection*{Local weak limit}
Let $\mathcal{P}(\mathcal{G}_{\star})$
denote the set of Borel
probability measures on~$\mathcal{G}_{\star}$, equipped with the
topology of weak
convergence \cite{Bil99}. Given a finite graph $G=(V,E)$, let
$\mathcal{U}(G)$
denote the law induced on $\mathcal{G}_{\star}$ by rooting $G$ at a
uniformly chosen
vertex $o\in V$ and restricting $G$ to the
connected component of $o$. If $\{G_n\}_{n\geq1}$ is a sequence of
finite graphs
such that $\{\mathcal{U}(G_n)\}_{n\geq1}$ admits a limit $\mu\in
\mathcal{P}(\mathcal{G}_{\star})$, we
call $\mu$ the \textit{local weak limit} of $\{G_n\}_{n\geq1}$ and write
\[
G_n \mathop{\longrightarrow}_{n\to\infty}^{\mathrm{LWC}} \mu.
\]
%

\subsection*{Edge-rooted graphs} Let $\mathcal{G}_{\star\star}$
denote the set of locally
finite connected graphs with a distinguished oriented edge, taken up to
the natural isomorphism relation and equipped with the natural
distance. With any function $f\dvtx\mathcal{G}_{\star\star}\to
\mathbb{R}$ is naturally
associated a function $\partial f \dvtx\mathcal{G}_{\star}\to
\mathbb{R}$, defined by
\[
\partial f(G,o)  =  \sum_{i\sim o}f(G,i,o).
\]
Dually, with any $\mu\in\mathcal{P}(\mathcal{G}_{\star})$ is
naturally associated a
nonnegative measure $\vec{\mu}$ on $\mathcal{G}_{\star\star}$,
defined by the following
relation: for any Borel $f\dvtx\mathcal{G}_{\star\star}\to
[0,\infty)$,
\begin{eqnarray*}
\int_{\mathcal{G}_{\star\star}}f\,d\vec{\mu}& = &\int_{\mathcal
{G}_{\star}}(
\partial f) \,d\mu.
\end{eqnarray*}
Note that $\vec{\mu}(\mathcal{G}_{\star\star})=\operatorname{deg} (\mu)$, where
$\operatorname{deg} (\mu)$ is the average
degree of the root:
\begin{eqnarray*}
\operatorname{deg}(\mu) & := & \int_{\mathcal{G}_{\star}}\operatorname{deg} (G,o)\,d\mu(G,o).
\end{eqnarray*}
\subsection*{Unimodularity}
Given $f\dvtx\mathcal{G}_{\star\star}\to\mathbb{R}$, we define
its \textit{reversal} $f^*\dvtx\mathcal{G}_{\star\star}
\to\mathbb{R}$ by
\begin{eqnarray*}
f^*(G,i,o) & = & f(G,o,i).
\end{eqnarray*}
It was shown in \cite{AldLyo07} that any $\mu\in\mathcal
{P}(\mathcal{G}_{\star})$ arising as
the local weak limit of some sequence of finite graphs satisfies the symmetry
\begin{eqnarray}
\label{eqmtp} \int_{\mathcal{G}_{\star\star}}f\,d\vec{\mu}& = & \int
_{\mathcal
{G}_{\star\star}} f^*\,d\vec{\mu},
\end{eqnarray}
for any Borel $f\dvtx\mathcal{G}_{\star\star}\to[0,\infty)$. A
measure $\mu\in\mathcal{P}(\mathcal{G}_{\star}
)$ satisfying (\ref{eqmtp}) is called \textit{unimodular}, and the set
of such measures is denoted by $\mathcal{U}$. The property (\ref
{eqmtp}) may
be viewed as an infinite analogue of the trivial identity
\begin{eqnarray*}
\sum_{o\in V}\sum_{i\sim o}f(i,o)
& = & \sum_{o\in V}\sum_{i\sim o}f(o,i),
\end{eqnarray*}
valid for any finite graph $G=(V,E)$ and any $f\dvtx\vec{E}\to
\mathbb{R}$.

\subsection*{Marks on oriented edges}
It will sometimes be convenient to work with \textit{networks}, that
is, graphs equipped with a map from $\vec{E}$ to some fixed complete
separable metric space $\Xi$. 
The above definitions extend naturally; see \cite{AldLyo07}.

\subsection*{Unimodular Galton--Watson trees}
Let $\pi=\{\pi_k\}_{k\geq0}$ be a probability distribution on
$\mathbb{N}$
with finite, nonzero mean.
A \textit{unimodular Galton--Watson tree} with degree distribution $\pi$
is a random rooted tree obtained by a Galton--Watson branching process
where the root has offspring distribution $\pi$ and all descendants
have the size-biased offspring distribution $\widehat{\pi}=\{\widehat
{\pi}_k\}_{k\geq0}$, where
\begin{eqnarray}
\label{eqsizebiased}
\widehat\pi_k & = & \frac{(k+1)\pi_{k+1}}{\sum_{i}i\pi_{i}}.
\end{eqnarray}
The resulting law is unimodular and is denoted by $\operatorname{UGWT}(\pi)$.
Such trees play a distinguished role in the local weak convergence
theory, as they are the limits of many natural sequences of random
graphs, including the Erd\H{o}s--R\'enyi model and more generally,
random graphs with
prescribed degrees.

\subsection*{The pairing model}
Given a sequence ${\mathbf d}=\{d(i)\}
_{1\leq
i\leq n}$ of nonnegative integers whose sum is even, the pairing model
\cite{Bol80,Jan09} generates a random graph $\mathbb{G}[\mathbf{d}]$
on $V=\{
1,\ldots,n\}$ as follows: $d(i)$ \textit{half-edges} are attached to each
$i\in V$, and the $2m=d(1)+\cdots+d(n)$ half-edges are paired uniformly
at random to form $m$ edges. Loops and multiple edges are then removed
(a few variants exist---see \cite{Hof13}---but they are equivalent
for our purpose). Now, consider a degree sequence ${\mathbf d}_n=\{
d_n(i)\}
_{1\leq i\leq n}$ for each $n\geq1$ and assume that
\begin{equation}
\label{h1}
\forall k\in\mathbb{N},\qquad\frac{1}{n}\sum
_{i=1}^n\mathbf{1}_{\{
d_n(i)=k\}
}\mathop{\longrightarrow}_{n\to\infty} \pi_k,
\end{equation}
for some probability distribution $\pi=\{\pi_k\}_{k\geq0}$ on
$\mathbb{N}$
with finite, nonzero mean. Under the additional assumption that
\[
\sup_{n\geq1} \Biggl\{\frac{1}{n}\sum
_{i=1}^n d^2_n(i) \Biggr\} <
\infty,
\]
the local weak limit of $\{\mathbb{G}[{\mathbf d_n}]\}_{n\geq1}$ is
$\mu:=\operatorname{UGWT}(\pi)$ almost surely; see \cite{Bor12}.

\section{Results}
Our first main result is that the notion of balanced allocations can be
extended from finite graphs to their local weak limits, in such a way
that the induced loads behave continuously with respect to local weak
convergence. Let us define a \textit{Borel allocation} as a measurable
function $\Theta\dvtx\mathcal{G}_{\star\star}\to[0,1]$ such that
$\Theta+\Theta^*=1$, and
call it \textit{balanced} on a given $\mu\in\mathcal{U}$ if for $\vec{\mu}$-almost-every $(G,i,o)\in\mathcal{G}_{\star\star}$,
\[
\partial\Theta(G,i)<\partial\Theta(G,o) \quad \Longrightarrow\quad \Theta (G,i,o)=0.
\]
This definition is the natural analogue of (\ref{eqbalance}) when
finite graphs are replaced by unimodular measures. We then have the
following result.

\begin{theorem}
\label{thmain1}
Let $\mu\in\mathcal{U}$ be such that $\operatorname{deg} (\mu)<\infty$. Then,
\begin{longlist}[{}]
\item[\textup{Existence and optimality}.] There is a Borel allocation $\Theta_0$
that is balanced on $\mu$, and for any Borel allocation $\Theta$ the
following are equivalent:
\begin{enumerate}[(iii)]
\item[(i)] $\Theta$ is balanced on $\mu$.
\item[(ii)] $\Theta$ minimizes $\int f\circ\partial\Theta\,d\mu$ for some
strictly convex $f\dvtx[0,\infty)\to\mathbb{R}$.
\item[(iii)] $\Theta$ minimizes $\int f\circ\partial\Theta\,d\mu$ for every
convex $f\dvtx[0,\infty)\to\mathbb{R}$.
\item[(iv)] $\partial\Theta=\partial\Theta_0$, $\mu$-almost-everywhere.
\end{enumerate}
\item[\textup{Continuity}.] For any sequence $\{G_n\}_{n\geq1}$ of finite graphs,
\[
 \bigl(G_n \mathop{\longrightarrow}_{n\to\infty}^{\mathrm{LWC}}
\mu \bigr) \quad \Longrightarrow\quad \bigl(\mathcal{L}_{G_n}\mathop{
\longrightarrow}_{n\to
\infty}^{\mathcal{P}(\mathbb{R})}\mathcal{L} _{\mu} \bigr),
\]
where $\mathcal{L}_\mu$ is the law of the random variable $\partial
\Theta_0\in
L^1(\mathcal{G}_{\star},\mu)$.
\item[\textup{Variational characterization}.]  The stop-loss transform of
$\mathcal{L}_\mu
$ is given by
\[
\Psi_{\mathcal{L}_\mu}(t)  =  \mathop{\max_{f\dvtx\mathcal{G}_{\star}\to[0,1]}}_{\mathrm{Borel}}
 \biggl\{
\frac{1}{2}\int_{\mathcal{G}_{\star\star}}\widehat{f}\, d\vec{\mu}- t\int
_{\mathcal{G}_{\star}}f\, d\mu \biggr\}, \qquad t\geq0,
\]
where $\widehat{f}(G,i,o):=f(G,o)\wedge f(G,i)$.
\end{longlist}
\end{theorem}

Recall that the \textit{stop-loss transform} of a nonnegative integrable
random variable $X$ (in fact, of its law $\mathcal{L}$) is the
function $\Psi
_X=\Psi_\mathcal{L}$ defined by
\begin{eqnarray*}
\Psi_X (t) & = & \mathbb{E}\bigl[(X-t)^+\bigr] = \int
_{\mathbb
{R}}(x-t)^+\,d\mathcal{L}(x), \qquad t\geq0.
\end{eqnarray*}
This function plays a central role in the theory of \textit{convex
ordering}, due to the classical equivalence between the following
conditions (see, e.g., \cite{ShaSha94}):
\begin{longlist}[(ii)]
\item[(i)] $\mathbb{E}[f(X)]\leq\mathbb{E}[f(Y)]$ for every convex
function $f\dvtx
[0,\infty)\to\mathbb{R}$.
\item[(ii)] $\Psi_X\leq\Psi_Y$ and $\Psi_X(0)=\Psi_Y(0)$.
\end{longlist}
In particular, if $\Psi_X=\Psi_Y$ then $X$ and $Y$ have the same
Laplace transforms. This shows that $\Psi_\mathcal{L}$ characterizes
$\mathcal{L}$.
Consequently, the above variational problem completely determines the
limiting empirical load distribution.

Our second main result is an explicit resolution of this variational
problem in the important case where $\mu=\operatorname{UGWT}(\pi)$, for an
arbitrary degree distribution $\pi=\{\pi_k\}_{k\geq0}$ on $\mathbb
{N}$ with
finite, nonzero mean. Throughout the paper, we let $[x]^1_0$ denote
the closest point to $x\in\mathbb{R}$ in the interval $[0,1]$, that is,
\begin{eqnarray*}
[x]^1_0 & := & %
\cases{ 0, & \quad$\mbox{if
}x \leq0$,\vspace*{2pt}
\cr
x, & \quad$\mbox{if }x\in[0,1]$,\vspace*{2pt}
\cr
1, &
\quad$\mbox{if }x \geq1$.}
\end{eqnarray*}
Given $t\in\mathbb{R}$ and $Q\in\mathcal{P}([0,1])$, we let $F_{\pi
,t}(Q)\in\mathcal{P}([0,1])$
denote the law of
\[
[1-t+\xi_1+\cdots+\xi_{\widehat{D}} ]^1_0,
\]
where ${\widehat{D}}$ follows the size-biased distribution $\widehat{\pi
}$ defined at (\ref{eqsizebiased}), and where $\{\xi_k\}_{k\geq1}$
are {i.i.d.} with law $Q$, independent of ${\widehat{D}}$. As
conjectured by Hajek \cite{Haj90}, the value of $\Psi_{\mathcal
{L}_\mu}(t)$
turns out to be controlled by the solutions to the distributional fixed
point equation $Q=F_{\pi,t}(Q)$. The latter can be solved numerically;
see \cite{Haj90} for detailed tables in the case where $\pi$ is Poisson.

\begin{theorem}
\label{thmain2}
When $\mu=\operatorname{UGWT}(\pi)$, we have for every $t\in\mathbb{R}$:
\begin{eqnarray*}
\Psi_{\mathcal{L}_\mu}(t)
& = & \max
_{Q=F_{\pi,t}(Q)} \biggl\{\frac{\mathbb{E}[D]}{2}\mathbb {P} (\xi_1+
\xi _2>1 )-t\mathbb{P} (\xi_1+\cdots+\xi_D>t
) \biggr\},
\end{eqnarray*}
where $D\sim\pi$ and where $\{\xi_k\}_{k\geq1}$ are {i.i.d.} with
law $Q$, independent of $D$. The maximum is over all choices of $Q\in
\mathcal{P}
([0,1])$ subject to $Q=F_{\pi,t}(Q)$.
\end{theorem}

By analogy with the case of finite graphs, we define the maximum
subgraph density of a measure $\mu\in\mathcal{U}$ with $\operatorname{deg} (\mu
)<\infty$ as
the essential supremum of the random variable $\partial\Theta_0$
constructed in Theorem~\ref{thmain1}. In other words,
\[
\varrho(\mu):=\sup\bigl\{t\in\mathbb{R}\dvtx\Psi_{\mathcal{L}_\mu
}(t)>0\bigr\}.
\]
In light of Theorem~\ref{thmain1}, it is natural to seek a continuity
principle of the form
\begin{equation}
\label{eqideal} \bigl(G_n \mathop{\longrightarrow}_{n\to\infty}^{\mathrm{LWC}}
\mu \bigr) \quad \Longrightarrow \quad \bigl(\varrho(G_n)\mathop{\longrightarrow}_{n\to
\infty}\varrho (\mu) \bigr).
\end{equation}
However, a moment of thought shows that the graph parameter $\varrho
(G)$ is too sensitive to be captured by local weak convergence. Indeed,
if $|V(G_n)|\to\infty$, then adding a large but fixed clique to $G_n$
will arbitrarily boost the value of $\varrho(G_n)$ without affecting
the local weak limit of $\{G_n\}_{n\geq1}$.
Similarly, for random graphs with a prescribed asymptotic degree
distribution $\pi\in\mathcal{P}(\mathbb{N})$, we expect (\ref
{eqideal}) to fail when
$\pi$ has heavy tail, due to the presence of extremely dense subgraphs
with negligible size. Understanding the maximum subgraph density in
that regime remains an interesting open question.
Our third main result establishes~(\ref{eqideal}) in the light-tail regime.

\begin{theorem}
\label{thmain3}
Consider a sequence $\{{\mathbf d_n}\}_{n\geq1}$ of degree sequences that
approach some distribution $\pi=\{\pi_k\}_{k\geq0}$ in the sense of
(\ref{h1}). Assume that $\pi_0+ \pi_1<1$ and that $\pi$ has light tail,
that is, for some $\theta>0$,
\begin{equation}
\label{h2}
\sup_{n\geq1} \Biggl\{\frac{1}{n}\sum
_{i=1}^ne^{\theta d_n(i)} \Biggr\} <\infty.
\end{equation}
Then $\varrho (\mathbb{G}[{\mathbf d_n}] )
{\ds\mathop{\longrightarrow}_{n\to\infty}^{\mathbb{P}}}\varrho(\mu)$,
where $\mu=\operatorname{UGWT}(\pi)$.
\end{theorem}

In particular, the result applies to the Erd\H{o}s--R\'enyi random
graph $\mathbb{G}_n$ with
$n$ vertices and $m=\lfloor\alpha n\rfloor$ edges. Indeed, the
conditional law of $\mathbb{G}_n$ given its (random) degree sequence
$\mathbf
{d_n}$ is precisely that of $\mathbb{G}[\mathbf{d_n}]$, and $\{
{\mathbf{d_n}}\}_{n\geq
1}$ satisfies a.s. the conditions (\ref{h1}) and (\ref{h2}) with $\pi
=\operatorname{Poisson}(2\alpha)$. Therefore,
Theorems \ref{thmain1}, \ref{thmain2} and \ref{thmain3} settle the
conjectures of \cite{Haj96} and validate the numerical tables for
$\varrho$ given therein. We note that the quantity $\varrho$ depends
monotonically and continuously on the connectivity parameter $\alpha$.
More precisely, it is not hard to show that for any positive $\alpha
<\beta$,
\[
1 \leq\frac{\varrho(\beta)}{\varrho(\alpha)} \leq\frac{\beta
}{\alpha}.
\]

Also, since a graph $G$ is \textit{$k$-orientable} ($k\in\mathbb{N}$)
if and only
if $\varrho(G)< k$, Theorem~\ref{thmain3} extends recent results on
the $k$-orientability of the Erd\H{o}s--R\'enyi random graph \cite{FerRam07,CaiSan07}.
See \cite{GaoWor10,FouKho13,Lel12,LecLel12} for various generalizations.

\section{Proof ingredients and related work}

\subsection*{The objective method}
This work is a new illustration of the general principles exposed in
the \textit{objective method} by Aldous and Steele \cite{AldSte04}. The
latter provides a powerful framework for the unified study of sparse
random graphs and has already led to several remarkable results. Two
prototypical examples are the celebrated $\zeta(2)$ limit in the random
assignment problem due to Aldous \cite{Ald01}, and the asymptotic
enumeration of spanning trees in large graphs by Lyons \cite{Lyo05}.
Since then, the method has been successfully applied to various other
combinatorial enumeration/optimization problems on graphs, including
(but not limited to) \cite{Ste02,GamNow06,SalSha09,BorLel13,Lel12,Sal13,LecLel12,KhaSun12}.

\subsection*{Lack of correlation decay}
In the problem considered here, there is a major obstacle to a direct
application of the objective method: the balanced load at a vertex is
\textit{not} determined by the local environment around that vertex. For
example, every vertex of a $d$-regular graph with girth $h$ has the
same $h$-neighborhood as the root of a $d$-regular tree with height
$h$. However, the induced load is $\frac{d}{2}$ in the first case and
$1-\frac{1}{(d-1)^{h-1}d}$ in the second. This long-range dependence
gives rise to nonuniqueness issues when trying to extend the notion of
balanced loads to infinite graphs. We refer to \cite{Haj96} for a
detailed study of this phenomenon---therein called \textit{load
percolation}---as well as several fascinating questions.

\subsection*{Relaxation} To overcome the lack of correlation decay, we
introduce a suitable relaxation of the balancing condition (\ref{eqbalance}), which we call $\varepsilon$-balancing. Remarkably
enough, any
positive value of the perturbative parameter $\varepsilon$ suffices to
annihilate the long-range dependences described above. This allows us
to define a unique $\varepsilon$-balanced Borel allocation $\Theta
_\varepsilon\dvtx\mathcal{G}_{\star\star}
\to[0,1]$ and to establish the continuity of the induced load
$\partial
\Theta_{\varepsilon}\dvtx\mathcal{G}_{\star}\to[0,\infty)$ with
respect to local convergence
(Section~\ref{seceps}). We then use unimodularity to prove that, as
the perturbative parameter $\varepsilon$ tends to $0$, $\Theta
_\varepsilon$ converges in
a certain sense to a balanced Borel allocation $\Theta_0$
(Section~\ref{seczero}). This quickly leads to a proof of Theorem~\ref
{thmain1}
(Section~\ref{secthmain1}). In spirit, the role of the perturbative
parameter $\varepsilon>0$ is comparable to that of the temperature in
\cite
{BorLel13}, although no Gibbs--Boltzmann measure is involved in the
present work. 

\subsection*{Recursion on trees} As many other graph-theoretical
problems, load balancing has a simple recursive structure when
considered on trees. Indeed, once the value of the allocation along a
given edge $\{i,j\}$ has been fixed, the problem naturally decomposes
into two independent sub-problems, corresponding to the two disjoint
subtrees formed by removing $\{i,j\}$. Note, however, that in the
resulting sub-problems the loads of $i$ and $j$ must be shifted by a
suitable amount to take into account the contribution of the removed
edge. The precise effect of this shift on the loads induced at $i$ and
$j$ defines what we call the \textit{response functions} of the two
subtrees (Section~\ref{secbaseload}). Those response functions satisfy
a recursion (Section~\ref{secrecursion}). Recursions on trees
automatically give rise to distributional fixed-point equations when
specialized to Galton--Watson trees. Such equations are a common
ingredient in the objective method; see \cite{AldBan05}. This leads to
the proof of Theorem~\ref{thmain2} (Section~\ref{secthmain2}).

\subsection*{Dense subgraphs in the pairing model} Finally, the proof of
Theorem~\ref{thmain3} (Section~\ref{secthmain3}) relies on a property
of random graphs with a prescribed degree sequence that might be of
independent interest: under the exponential moment assumption~(\ref
{h2}), we show that dense subgraphs must be extensively large with high
probability. Our argument is based on the first-moment method. See
Proposition~\ref{prpairing} for the precise statement, and \cite{Luc92}, Lemma~6,  for a result in the same direction.

\section{\texorpdfstring{$\varepsilon$}{varepsilon}-balancing}
\label{seceps}Throughout this section, $G=(V,E)$ is a locally finite
graph and $\varepsilon>0$ is a fixed parameter. An allocation $\theta$
on $G$ is called \textit{$\varepsilon$-balanced} if for every $(i,j)\in
\vec{E}$,
\begin{eqnarray}
\label{eqepsbb}
\theta(i,j) & = & \biggl[\frac{1}{2}+\frac{\partial\theta
(i)-\partial
\theta(j)}{2\varepsilon}
\biggr]^1_0.
\end{eqnarray}
This condition can be viewed as a relaxation of (\ref{eqbalance}). Its
interest lies in the fact that it fixes the nonuniqueness issue on
infinite graphs.

\begin{proposition}[(Existence, uniqueness and monotony)]
\label{preps}
If $G$ has\break bounded degrees, then there is a unique $\varepsilon$-balanced allocation $\theta$ on $G$. If moreover $E'\subseteq E$,
then the $\varepsilon$-balanced allocation $\theta'$ on $G'=(V,E')$
satisfies $\partial\theta'\leq\partial\theta$.
\end{proposition}

\begin{pf}
The set $K$ of all allocations on $G$ is clearly a compact convex
subset of the locally convex space $\mathbb{R}^{\vec{E}}$ equipped
with the
topology of coordinate-wise convergence. Moreover, the mapping $K\ni
\theta\mapsto\theta'\in K$ defined by
\begin{eqnarray*}
\theta'(i,j) & = & \biggl[\frac{1}{2}+
\frac{\partial\theta
(i)-\partial
\theta(j)}{2\varepsilon} \biggr]^1_0
\end{eqnarray*}
is continuous. It must therefore admit a fixed point, by the
Schauder--Tychonoff fixed-point theorem (see, e.g., \cite{AgaMee01}, Theorem~8.2). This proves existence. Now, consider
$E'\subseteq E$
and let $\theta,\theta'$ be $\varepsilon$-balanced allocations on
$G,G'$, respectively. Fix $o\in V$ and set
\[
I:=\bigl\{i\in V\dvtx\{i,o\}\in E', \theta'(i,o)>
\theta(i,o)\bigr\}.
\]
Clearly,
\begin{eqnarray*}
\partial\theta'(o)-\partial\theta(o) & \leq& \sum
_{i\in I} \bigl(\theta '(i,o)-\theta(i,o)
\bigr).
\end{eqnarray*}
On the other hand, since the map $x\mapsto [\frac{1}{2}+\frac
{x}{2\varepsilon} ]^1_0$ is nondecreasing and $\frac
{1}{2\varepsilon}$-Lipschitz, our assumption on $\theta,\theta'$
implies that for every $i\in I$,
\begin{eqnarray*}
\theta'(i,o)-\theta(i,o) & \leq& \frac{1}{2\varepsilon} \bigl(\partial
\theta'(i)-\partial\theta(i)-\partial\theta'(o)+\partial
\theta(o) \bigr).
\end{eqnarray*}
Injecting this into the above inequality and rearranging, we obtain
\begin{eqnarray}
\nonumber
\partial\theta'(o)-\partial\theta(o) & \leq&
\frac{1}{|I|+
2\varepsilon
}\sum_{i\in I} \bigl(\partial
\theta'(i)-\partial\theta(i) \bigr)
\nonumber
\\[-8pt]
\label{eqcontraction} \\[-8pt]
\nonumber
& \leq& \frac{\Delta}{\Delta+ 2\varepsilon}\max
_{i\in I} \bigl(\partial \theta'(i)-\partial\theta(i)
\bigr),
\end{eqnarray}
where $\Delta$ denotes the maximum degree in $G$. Now, observe that
$\partial\theta,\partial\theta'$ are $[0,\Delta]$-valued, so that
$M:=\sup_{V}(\partial\theta'-\partial\theta)$ is finite. Property
(\ref{eqcontraction}) forces $M\leq0$, which proves the monotony
$E'\subseteq E\Longrightarrow\partial\theta'\leq\partial\theta$. In
particular, $E'=E$ implies $\partial\theta' = \partial\theta$,
which in
turns forces $\theta'=\theta$, thanks to (\ref{eqcontraction}).
\end{pf}
We now remove the bounded-degree assumption as follows. Fix $\Delta\in
\mathbb{N}$, and consider the truncated graph $G^\Delta=(V,E^\Delta)$ obtained
from $G$ by isolating all nodes having degree more than $\Delta$, that is,
\[
E^\Delta= \bigl\{\{i,j\}\in E\dvtx\operatorname{deg} (G,i)\vee\operatorname{deg} (G,j)\leq \Delta
\bigr\}.
\]
By construction, $G^\Delta$ has degree at most $\Delta$, and we let
$\Theta^\Delta_\varepsilon(G,i,j)$ denote the mass sent along
$(i,j)\in
\vec{E}$ in the unique $\varepsilon$-balanced allocation on $G^\Delta
$, with the
understanding that $\Theta^\Delta_\varepsilon(G,i,j)=0$ if $\{i,j\}
\notin
E^\Delta$.
By uniqueness, this quantity depends only on the isomorphism class of
the edge-rooted graph $(G,i,j)$, so that we have a well-defined map
$\Theta^\Delta_{\varepsilon}\dvtx\mathcal{G}_{\star\star}\to
[0,1]$. By an immediate
induction on $r\in\mathbb{N}$, the local contraction (\ref
{eqcontraction}) yields
\[
[G,o]_r \equiv\bigl[G',o'
\bigr]_r \quad \Longrightarrow\quad \bigl\llvert \partial\Theta
^\Delta _\varepsilon(G,o)-\partial\Theta^\Delta_\varepsilon
\bigl(G',o'\bigr)\bigr\rrvert \leq \Delta \biggl(1+
\frac{2\varepsilon}{\Delta} \biggr)^{-r}.
\]
Since the map $x\mapsto [\frac{1}2+\frac{x}{2\varepsilon}
]^1_0$ is $\frac{1}{2\varepsilon}$-Lipschitz, it follows that
\[
[G,i,j]_r \equiv\bigl[G',i',j'
\bigr]_r \quad \Longrightarrow\quad \bigl\llvert \Theta^\Delta_\varepsilon(G,i,j)-
\Theta ^\Delta _\varepsilon\bigl(G',i',j'
\bigr)\bigr\rrvert \leq\frac{\Delta}{2\varepsilon
} \biggl(1+\frac{2\varepsilon}{\Delta}
\biggr)^{-r}\!.
\]
Thus, the map $\Theta^\Delta_\varepsilon$ is equicontinuous.
Now, the sequence of sets $\{E_\Delta\}_{\Delta\geq1}$ increases to
$E$, so the monotony in Proposition~\ref{preps} guarantees that $\{
\partial\Theta^\Delta_\varepsilon\}_{\Delta\geq1}$ converges
pointwise on $\mathcal{G}_{\star}$. Moreover, any given $\{i,j\}\in E$ belongs to $E^\Delta$ for large
enough $\Delta$, and the definition of $\varepsilon$-balancing yields
\[
\Theta^\Delta_\varepsilon(G,i,j)  =  \biggl[
\frac{1}{2}+\frac
{\partial\Theta
^\Delta_\varepsilon(G,i)-\partial\Theta^\Delta_\varepsilon
(G,j)}{2\varepsilon} \biggr]^1_0.
\]
Consequently, the pointwise limit $\Theta_\varepsilon:=\lim_{\Delta
\to\infty
}\Theta^\Delta_\varepsilon$ exists in $[0,1]^{\mathcal{G}_{\star
\star}}$. It clearly satisfies
$\Theta_\varepsilon+\Theta^*_\varepsilon=1$ and it is Borel as the
pointwise limit of
continuous maps. Thus, it is a Borel allocation. Letting $\Delta\to
\infty$ above yields
\begin{equation}
\label{eqepsbalance}
\Theta_\varepsilon(G,i,j)  =  \biggl[\frac{1}{2}+
\frac{\partial
\Theta_\varepsilon
(G,i)-\partial\Theta_\varepsilon(G,j)}{2\varepsilon} \biggr]^1_0.
\end{equation}

\section{The \texorpdfstring{$\varepsilon\to0$}{varepsilonto0} limit}
\label{seczero}
We now send the perturbative parameter $\varepsilon$ to $0$, and show that
$\Theta_\varepsilon$ converges in a certain sense to a balanced Borel
allocation
$\Theta_0$.
Fix $\mu\in\mathcal{U}$ with $\operatorname{deg} (\mu)<\infty$. We write $\|f\|
_p$ for the
norm in both $L^p(\mu)$ and $L^p(\vec{\mu})$: which is meant should
be clear
from the context. Note that by unimodularity, we have for any Borel
allocation $\Theta$,
\begin{equation}
\label{ineqL1}
\|\Theta\|_{1}  =  \int_{\mathcal{G}_{\star\star}}
\Theta\, d\vec{\mu} =  \int_{\mathcal{G}_{\star\star}}\frac
{\Theta+\Theta^*}{2}\, d\vec{
\mu} =  \frac{\operatorname{deg} (\mu)}{2}.
\end{equation}

\begin{proposition}[(Existence of a balanced Borel allocation)]
The limit ${\Theta_0}:=\lim_{\varepsilon\to0}\Theta_{\varepsilon}$
exists in $L^2(\vec{\mu})$ and is a balanced Borel allocation on $\mu$.
\label{przero}
\end{proposition}

\begin{pf}
We will establish the following property: for $0<\varepsilon\leq
\varepsilon'$,
\begin{equation}
\label{eqcauchy}
\|\Theta_{\varepsilon'}-\Theta_{\varepsilon}
\|^2_{2}  \leq \| \Theta _\varepsilon
\|^2_{2}-\|\Theta_{\varepsilon'}\|^2_{2}.
\end{equation}
In particular, $\|\Theta_\varepsilon\|^2_{2}\geq\|\Theta
_{\varepsilon
'}\|^2_{2}$ so $\lim_{\varepsilon\to0}\uparrow\|\Theta_\varepsilon
\|
^2_{2}$ exists. Consequently, the right-hand side tends to $0$ as
$\varepsilon,\varepsilon
'\to0$, hence so does the left-hand side. This provides a Cauchy
criterion in $L^2(\vec{\mu})$ for $\{\Theta_\varepsilon\}
_{\varepsilon>0}$, ensuring
the existence of $\Theta_0=\lim_{\varepsilon\to0}\Theta
_{\varepsilon
}$. The rest of the claim follows, since Borel allocations are closed
in $L^2(\vec{\mu})$ and letting $\varepsilon\to0$ in (\ref
{eqepsbalance}) shows
that $\Theta_0$ is balanced on $\mu$. In order to prove (\ref
{eqcauchy}), let us first assume that
\begin{equation}
\label{eqextra}
\mu \bigl(\bigl\{(G,o)\dvtx\operatorname{deg} (G,o)\leq\Delta\bigr\} \bigr)=1,
\end{equation}
for some $\Delta\in\mathbb{N}$.
This ensures that $f\in L^2(\vec{\mu})$, where
\[
f(G,i,o):=\partial\Theta_\varepsilon(G,o)+\varepsilon\Theta
_\varepsilon(G,i,o).
\]
A straightforward manipulation of (\ref{eqepsbalance}) shows that
\[
f(G,i,o)>f(G,o,i) \quad \Longrightarrow\quad \Theta_\varepsilon(G,i,o)=0.
\]
This implies $\Theta_\varepsilon f+\Theta^*_\varepsilon f^* = f\wedge
f^*$. On the other
hand, $f\wedge f^*\leq\Theta_{\varepsilon'} f+\Theta^*_{\varepsilon
'}f^*$ since $\Theta
_{\varepsilon'}+\Theta^*_{\varepsilon'}=1$. Thus,
$\Theta_\varepsilon f+\Theta^*_\varepsilon f^*
\leq
\Theta_{\varepsilon'} f+\Theta^*_{\varepsilon'} f^*$.
Integrating against $\vec{\mu}$ and invoking unimodularity, we get
$\langle
\Theta_\varepsilon-\Theta_{\varepsilon'},f\rangle_{L^2(\vec{\mu
})}\leq0$ or more explicitly,
\[
\langle\partial\Theta_{\varepsilon}-\partial\Theta_{\varepsilon
'},\partial
\Theta _{\varepsilon}\rangle_{L^2(\mu)} + \varepsilon\langle\Theta
_{\varepsilon
}-\Theta_{\varepsilon'}, \Theta_\varepsilon\rangle_{L^2(\vec{\mu
})}
 \leq 0.
\]
But we have not yet used $\varepsilon\leq\varepsilon'$, so we may
exchange $\varepsilon,\varepsilon'$ to get
\[
\langle\partial\Theta_{\varepsilon'}-\partial\Theta_{\varepsilon
},\partial
\Theta_{\varepsilon
'}\rangle_{L^2(\mu)} + \varepsilon' \langle
\Theta_{\varepsilon'}- \Theta_{\varepsilon}, \Theta_{\varepsilon'}
\rangle_{L^2(\vec{\mu})}  \leq 0.
\]
Adding-up those inequalities and rearranging, we finally arrive at
\[
\bigl(\varepsilon'-\varepsilon\bigr)\langle
\Theta_{\varepsilon}- \Theta _{\varepsilon'}, \Theta_{\varepsilon
'}
\rangle_{L^2(\vec{\mu})}  \geq \|\partial\Theta_{\varepsilon
}-\partial\Theta
_{\varepsilon'}\|_2^2 + \varepsilon\|
\Theta_{\varepsilon}-\Theta _{\varepsilon'}\|_2^2.
\]
In particular, $\langle\Theta_{\varepsilon}, \Theta_{\varepsilon'}
\rangle_{L^2(\vec{\mu})}
\geq\|\Theta_{\varepsilon'}\|^2_2$ and (\ref{eqcauchy}) follows since
\[
\|\Theta_{\varepsilon'}-\Theta_{\varepsilon'}\|^2_2  =
 \|\Theta _{\varepsilon'}\|^2_2+\|\Theta _{\varepsilon}
\|^2_2-2\langle\Theta_{\varepsilon}, \Theta
_{\varepsilon'} \rangle_{L^2(\vec{\mu})}.
\]
Finally, if our extra assumption (\ref{eqextra}) is dropped, we may
apply (\ref{eqcauchy}) with $\Theta_\varepsilon$, $\Theta
_{\varepsilon'}$
replaced by $\Theta^\Delta_\varepsilon$, $\Theta^\Delta
_{\varepsilon'}$ and let
then $\Delta\to\infty$. By construction, $\Theta^\Delta
_\varepsilon\to
\Theta_\varepsilon$ and $\Theta^\Delta_{\varepsilon'}\to\Theta
_{\varepsilon'}$
pointwise, and (\ref{eqcauchy}) follows by dominated convergence.
\end{pf}

\begin{proposition}[(Continuity of balanced loads)]
\label{prcontinuity}
Let $\{G_n\}_{n\geq1}$ be a sequence of finite graphs with local weak
limit $\mu$. Then
\[
\mathcal{L}_{G_n}\mathop{\longrightarrow}_{n\to\infty}^{\mathcal
{P}(\mathbb{R})}
\mathcal{L},
\]
where $\mathcal{L}=\mathcal{L}_\mu$ is the law of the random
variable $\partial\Theta
_0\in L^1(\mu)$.
\end{proposition}

\begin{pf}
For $n\geq1$ we let
$\widehat{G_n}$ denote the network obtained by encoding a balanced
allocation $\theta_n$ as $[0,1]$-valued marks on the oriented edges of
$G_n$. The sequence $\{\mathcal{U}(\widehat{G_n})\}_{n\geq1}$ is tight,
because $\{\mathcal{U}(G_n)\}_{n\geq1}$ converges weakly and the
marks are
$[0,1]$-valued. Consider any subsequential weak limit $ (\mathbb{G}
,o,\theta )$. By construction, $(\mathbb{G},o)$ has law $\mu$
and $\theta
$ is a.s. a balanced allocation on $\mathbb{G}$. Our goal is to
establish that
$\partial\theta(o) = \partial\Theta_0(\mathbb{G},o)$ a.s.
Set $\theta'(i,j):=\Theta_0(\mathbb{G},i,j)$. The random rooted
network $(\mathbb{G}
,o,\theta,\theta')$ is unimodular, since $(\mathbb{G},o,\theta)$ is
a weak
limit of finite networks and $\Theta_0$ is Borel. Now,
\begin{eqnarray*}
\mathbb{E} \bigl[ \bigl(\partial\theta(o)-\partial\theta'(o)
\bigr)^+ \bigr] & = & \mathbb{E} \biggl[\sum_{i\sim o}
\bigl(\theta(i,o)-\theta '(i,o) \bigr){\mathbf{1} }_{\partial\theta(o)>\partial\theta'(o)}
\biggr]
\\
& = & \mathbb{E} \biggl[\sum_{i\sim o} \bigl(
\theta(o,i)-\theta '(o,i) \bigr){\mathbf{1} }_{\partial\theta(i)>\partial\theta'(i)} \biggr]
\\
& = & \mathbb{E} \biggl[\sum_{i\sim o} \bigl(
\theta'(i,o)-\theta (i,o) \bigr){\mathbf{1} }_{\partial\theta(i)>\partial\theta'(i)}
\biggr],
\end{eqnarray*}
where the second equality follows from unimodularity and the third one
from the identities $\theta(o,i)=1-\theta(i,o)$ and $\theta
'(o,i)=1-\theta'(i,o)$. Combining the first and last lines, we see that
$\mathbb{E} [ (\partial\theta(o)-\partial\theta
'(o) )^+ ]$ equals
\[
\frac{1}{2} \mathbb{E} \biggl[\sum_{i\sim o}
\bigl(\theta(i,o)-\theta'(i,o) \bigr) ({\mathbf{1}
}_{\partial\theta(o)>\partial\theta'(o)}-{\mathbf{1}}_{\partial
\theta
(i)>\partial\theta'(i)} ) \biggr].
\]
The fact that $\theta,\theta'$ are balanced across $\{i,o\}$ easily
implies that $\theta(i,o)-\theta'(i,o)$ and ${\mathbf{1}}_{\partial
\theta
(o)>\partial\theta'(o)}
-{\mathbf{1}}_{\partial\theta(i)>\partial\theta'(i)}$ can neither be
simultaneously positive, nor simultaneously negative. Therefore, $
\mathbb{E} [ (\partial\theta(o)-\partial\theta'(o)
)^+ ]
\leq0$.
Exchanging the roles of $\theta,\theta'$ yields $\partial\theta
(o)=\partial\theta'(o)$ a.s., as desired.
\end{pf}

\section{Proof of Theorem~\texorpdfstring{\protect\ref{thmain1}}{1}}\label{secthmain1}
We now complete the proof of Theorem~\ref{thmain1}.

\begin{proposition}
\label{prdual}
Let $\Theta$ be a Borel allocation. Then for all $t\in\mathbb{R}$,
\begin{eqnarray*}
\int_{\mathcal{G}_{\star}} (\partial\Theta-t )^+\, d\mu & \geq& \mathop{\sup_{f\dvtx\mathcal{G}_{\star}\to[0,1]}}_{\mathrm{Borel}}
\biggl\{\frac{1}{2}\int_{\mathcal{G}_{\star\star}}
\widehat{f}\, d\vec{\mu}- t\int_{\mathcal{G}_{\star}}f\, d\mu \biggr\},
\end{eqnarray*}
with equality for all $t\in\mathbb{R}$ if and only if $\Theta$ is
balanced on
$\mu$.
\end{proposition}
\begin{pf}
Fix a Borel $f\dvtx\mathcal{G}_{\star}\to[0,1]$. Since $(\partial
\Theta-t)^+\geq
(\partial\Theta-t)f$, we have
\begin{eqnarray}
\label{ineqdual1}
\int_{\mathcal{G}_{\star}} (\partial\Theta-t )^+\, d\mu &
\geq& \int_{\mathcal{G}_{\star}}f\partial\Theta\, d\mu- t\int
_{\mathcal{G}_{\star}}f\, d\mu.
\end{eqnarray}
Using the unimodularity of $\mu$ and the identity $\Theta+\Theta^*=1$,
we also have
\begin{eqnarray}
\nonumber
\qquad\int_{\mathcal{G}_{\star}}f\partial\Theta\, d\mu & = &
\frac{1}{2}\int_{\mathcal{G}_{\star\star}} \bigl(f(G,o)\Theta(G,i,o)+f(G,i)
\Theta (G,o,i) \bigr)\, d\vec{\mu}(G,i,o)\hspace*{-10pt}
\nonumber
\\[-8pt]
\label{ineqdual2}
\\[-8pt]
\nonumber
 & \geq&\frac{1}{2}\int_{\mathcal{G}_{\star\star}}
\bigl(f(G,o)\wedge f(G,i) \bigr)\, d\vec{\mu}(G,i,o).
\end{eqnarray}
Combining (\ref{ineqdual1}) and (\ref{ineqdual2}) yields the
inequality. Let us examine the equality case. First, equality holds in
(\ref{ineqdual1}) if and only if for $\mu$-a.e. $(G,o)\in\mathcal
{G}_{\star}$,
\begin{eqnarray*}
\partial\Theta(G,o)>t\quad  & \Longrightarrow& \quad f(G,o)=1,
\\
\partial\Theta(G,o)<t \quad & \Longrightarrow & \quad f(G,o)=0.
\end{eqnarray*}
Second, equality holds in (\ref{ineqdual2}) if and only if for $\vec
{\mu}$-a.e. $(G,i,o)\in\mathcal{G}_{\star\star}$,
\begin{eqnarray*}
f(G,i)<f(G,o)\quad & \Longrightarrow& \quad \Theta(G,i,o)=0.
\end{eqnarray*}
If $\Theta$ is balanced on $\mu$, then the choice $f={\mathbf{1}}_{\{
\partial
\Theta>t\}}$ clearly satisfies all those requirements, so that equality
holds for each $t\in\mathbb{R}$ in the proposition. This proves the
\textit{if}
part and shows that the supremum in Proposition~\ref{prdual} is
attained, because at least one balanced allocation exists by
Proposition~\ref{przero}. Now, for the \textit{only if} part, suppose
that equality is achieved in Proposition~\ref{prdual}. Then the above
requirements imply that for $\vec{\mu}$-a.e. $(G,i,o)\in\mathcal
{G}_{\star\star}$,
\[
\partial\Theta(G,i)<t<\partial\Theta(G,o)\quad  \Longrightarrow \quad\Theta(G,i,o)=0.
\]
Since this must be true for all $t\in\mathbb Q$, it follows that
$\Theta
$ is balanced on $\mu$.
\end{pf}

\begin{pf*}{Proof of Theorem~\protect\ref{thmain1}}
Existence, continuity and the variational characterization were
established in Propositions \ref{przero}, \ref{prcontinuity} and \ref
{prdual}, respectively. Now, let $\Theta,\Theta'$ be Borel
allocations, and assume that $\Theta$ is balanced.
Applying Proposition~\ref{prdual} to $\Theta$ and $\Theta'$ shows that
for all $t\in\mathbb{R}$,
\begin{eqnarray*}
\int_{\mathcal{G}_{\star}}(\partial\Theta-t)^+ \, d\mu& \leq& \int
_{\mathcal{G}_{\star}}\bigl(\partial \Theta'-t\bigr)^+ \, d\mu.
\end{eqnarray*}
Moreover, (\ref{ineqL1}) guarantees that $\partial\Theta,\partial
\Theta'$ have the same mean. As already mentioned below the statement of
Theorem~\ref{thmain1}, those two conditions imply
\begin{eqnarray*}
\int_{\mathcal{G}_{\star}} (f\circ\partial\Theta ) \, d\mu& \leq& \int
_{\mathcal{G}_{\star}} \bigl(f\circ\partial\Theta' \bigr) \, d\mu,
\end{eqnarray*}
for any convex function $f\dvtx[0,\infty)\to\mathbb{R}$. We have
just proved
(i) $\Longrightarrow$ (iii). On the other hand, (iii) $\Longrightarrow$
(ii) is obvious. In particular, $\Theta_0$ satisfies (ii) and (iii).
The \textit{only if} part of Proposition~\ref{prdual} shows that (iii)
$\Longrightarrow$ (i). The implication (iv) $\Longrightarrow$ (iii) is
obvious given that $\Theta_0$ satisfies (iii). Thus, it only remains to
prove (ii) $\Longrightarrow$ (iv). Assume that $\Theta$ minimizes
$\int
(f\circ\partial\Theta)\,d\mu$ for some strictly convex function
$f\dvtx
[0,\infty)\to\mathbb{R}$, and let $m$ denote the value of this
minimum. Since
$\Theta_0$ satisfies (ii), we also have $\int(f\circ\partial
\Theta
_0)\,d\mu=m$. But then $\Theta':=(\Theta_0+\Theta)/{2}$ is an
allocation and by convexity,
\begin{eqnarray*}
\int_{\mathcal{G}_{\star}} \bigl(f\circ\partial\Theta'\bigr)\, d
\mu& \leq & \int_{\mathcal{G}_{\star}}\frac
{(f\circ\partial\Theta)+(f\circ\partial\Theta_0)}{2}\, d\mu = m.
\end{eqnarray*}
This inequality contradicts the definition of $m$, unless it is an equality.
This forces $\partial\Theta=\partial\Theta_0$ $\mu$-a.e., since $f$ is
strictly convex.
\end{pf*}

\section{Response functions}
\label{secbaseload}
As many other graph-theoretical problems, load balancing has a simple
recursive structure when specialized to trees. However, the exact
formulation of this recursion requires the possibility to \textit{condition} the allocation to take a certain value at a given edge, and
we first need to give a proper meaning to this operation.
Let $G=(V,E)$ be a locally finite graph and $b\dvtx V\to\mathbb{R}$
a function
called the \textit{baseload}. An allocation $\theta$ is \textit{balanced
with respect to} $b$ if
\begin{eqnarray*}
b(i)+\partial\theta(i)<b(j)+\partial\theta(j)\quad  & \Longrightarrow&\quad
\theta(i,j)=0,
\end{eqnarray*}
for all $(i,j)\in\vec{E}$. This is precisely the definition of
balancing, except that the load \textit{felt} by each vertex $i\in V$ is
shifted by a certain amount $b(i)$. Similarly, $\theta$ is \textit{$\varepsilon$-balanced with respect to} $b$ if for all $(i,j)\in\vec{E}$,
\begin{eqnarray}
\label{eqbalb}
\theta(i,j) & = & \biggl[\frac{1}{2}+\frac{b(i)+\partial\theta
(i)-b(j)-\partial\theta(j)}{2\varepsilon}
\biggr]^1_0.
\end{eqnarray}
The arguments used in Proposition~\ref{preps} are easily extended to
this situation.

\begin{proposition}[(Existence, uniqueness and  monotony)]
\label{prbaseload}
If $G$ has\break bounded degree and if $b$ is bounded, then there is a unique
$\varepsilon$-balanced allocation with baseload $b$. Moreover, if
$b'\leq b$ is
bounded and if $E'\subseteq E$, then the $\varepsilon$-balanced allocation
$\theta'$ on $G'=(V,E')$ with baseload $b'$ satisfies $b'+\partial
\theta
'\leq b+\partial\theta$.
\end{proposition}

As in Section~\ref{seceps}, we then define an $\varepsilon$-balanced
allocation
in the general case by considering the truncated graph $G^\Delta$ with
baseload the truncation of $b$ to $[-\Delta,\Delta]$, and let then
$\Delta\to\infty$. Monotony guarantees the existence of a limiting
$\varepsilon
$-balanced allocation. We shall need the following property.

\begin{proposition}[(Nonexpansion)]
\label{prnonexp}
Let $\theta,\theta'$ be the $\varepsilon$-balanced allocations with baseloads
$b,b'\dvtx V\to\mathbb{R}$. Set $f=\partial\theta+b$ and
$f'=\partial\theta
'+b'$. Then
\begin{eqnarray*}
\bigl\|f'-f \bigr\|_{\ell^1(V)} & \leq& \bigl\|b'-b
\bigr\|_{\ell^1(V)}.
\end{eqnarray*}
\end{proposition}

\begin{pf}
By considering $b''=b\wedge b'$ and using the triangle inequality, we
may assume that $b\leq b'$. Note that this implies $f\leq f'$, thanks
to Proposition~\ref{prbaseload}. When $G$ is finite, the claim
trivially follows from conservation of mass:
\begin{eqnarray*}
\sum_{o\in V} \bigl(f'(o)-f(o)
\bigr) & = & \sum_{o\in V} \bigl(b'(o)-b(o)
\bigr). 
\end{eqnarray*}
This then extends to the case where $G$ has bounded degrees with $b,b'$
bounded as follows: choose finite subsets $V_1\subseteq V_2\subseteq
\cdots$ such that $\bigcup_{n\geq1}V_n=V$. For each $n\geq1$, let
$\theta
_n,\theta_n'$ denote the $\varepsilon$-balanced allocations on the
subgraph induced by $V_n$, with baseloads the restrictions of $b,b'$ to
$V_n$. Then $\theta_n\to\theta$ and $\theta_n'\to\theta'$
pointwise, by
compactness and uniqueness. Now, any finite $K\subseteq V$ is contained
in $V_n$ for large enough $n$, and since $V_n$ is finite we know that
$f_n:=\partial\theta_n+b$ and $f'_n:=\partial\theta'_n+b'$ satisfy
\begin{eqnarray*}
\sum_{i\in K}\bigl|f_n'(i)-f_n(i)\bigr|
& \leq& \sum_{i\in V_n}\bigl|b'(i)-b(i)\bigr|.
\end{eqnarray*}
Letting $n\to\infty$ yields the desired result, since $K$ is arbitrary.
Finally, for the general case, we may apply the result to the truncated
graph $G^\Delta$ with baseloads the truncation of $b,b'$ to $[-\Delta
,\Delta]$, and let then $\Delta\to\infty$.
\end{pf}

Although the uniqueness in Proposition~\ref{prbaseload} does not
extend to the $\varepsilon=0$ case, the following weaker result will
be useful
in the next section.

\begin{proposition}[(Weak uniqueness)]\label{lmdicho}Assume that $\theta
,\theta'$ are balanced with respect to $b$ and that $\|\partial\theta
-\partial\theta'\|_{\ell^1(V)}<\infty$. Then, $\partial\theta
=\partial
\theta'$.
\end{proposition}

\begin{pf}
Fix $\delta>0$. Since $\|\partial\theta-\partial\theta'\|_{\ell
^1(V)}<\infty$, the level set
$S:=\{j\in V\dvtx\partial\theta'(j)-\partial\theta(j)>\delta\}$
must be finite. Therefore, it satisfies the conservation of mass:
\begin{eqnarray}
\label{eq03}
\sum_{j\in S}\partial
\theta'(j)-\partial\theta(j) & = & \sum_{(i,j)\in
E(V-S,S)}
\theta'(i,j)-\theta(i,j). 
\end{eqnarray}
Now, if $(i,j)\in E(V-S,S)$ then clearly, $\partial\theta
'(i)-\partial
\theta(i)< \partial\theta'(j)-\partial\theta(j)$.
Consequently, at least one of the following inequalities must hold:
\[
b(j)-b(i)<\partial\theta(i)-\partial\theta(j)\quad\mbox{or}\quad
b(j)-b(i)> \partial\theta'(i)-\partial\theta'(j).
\]
The first one implies $\theta(i,j)=1$ and the second $\theta'(i,j)=0$,
since $\theta,\theta'$ are balanced with respect to $b$. In either
case, we have $\theta'(i,j)\leq\theta(i,j)$. Thus, the right-hand side
of (\ref{eq03}) is nonpositive, hence so must the left-hand side be.
This contradicts the definition of $S$ unless $S=\varnothing$, that is,
$\partial\theta'\leq\partial\theta+\delta$. Since $\delta$ is
arbitrary, we conclude that $\partial\theta'\leq\partial\theta$.
Equality follows by symmetry.
\end{pf}

Given $o\in V$ and $x\in\mathbb{R}$, we set $\mathfrak
{f}^\varepsilon_{(G,o)}(x)=x+\partial\theta
(o)$ where $\theta$ is the $\varepsilon$-balanced allocation with
baseload $x$
at $o$ and $0$ elsewhere.
We call $\mathfrak{f}^\varepsilon_{(G,o)}\dvtx\mathbb{R}\to
\mathbb{R}$ the \textit{response function} of
the rooted graph $(G,o)$. Propositions \ref{prbaseload} and \ref
{prnonexp} guarantee that $\mathfrak{f}^\varepsilon_{(G,o)}$ is
nondecreasing and
nonexpansive, that is,
\begin{eqnarray}
\label{eqnondnone} x\leq y \quad  & \Longrightarrow & \quad 0\leq\mathfrak{f}^\varepsilon
_{(G,o)}(y)-\mathfrak{f}^\varepsilon_{(G,o)}(x)\leq y-x.
\end{eqnarray}
Note for future use that the definition of $\mathfrak{f}^\varepsilon
_{(G,o)}(x)$ also implies
\begin{eqnarray}
\label{eqcoercive} 0 \leq& \mathfrak{f}^\varepsilon_{(G,o)}(x)-x & \leq
\operatorname{deg} (G,o).
\end{eqnarray}
When $G$ is a tree, response functions turn out to satisfy a simple recursion.

\section{Recursion on trees}
\label{secrecursion}
We are now ready to state the promised recursion. Fix a tree $T=(V,E)$.
Deleting $\{i,j\}\in{E}$ creates two disjoint subtrees, viewed as
rooted at $i$ and $j$ and denoted $T_{i\to j}$ and $T_{j\to i}$, respectively.

\begin{proposition}
\label{prreceps}
The response function $\mathfrak{f}^{\varepsilon}_{(T,o)}$ is
invertible and
\begin{eqnarray}
\label{eqreceps} \bigl\{\mathfrak{f}^{\varepsilon}_{(T,o)} \bigr
\}^{-1} & = & \mathrm{Id}-\sum_{i\sim
o} \bigl[1-
\bigl\{\mathfrak{f}^{\varepsilon}_{T_{i\to
o}}+\varepsilon(2\mathrm{Id}-1) \bigr
\}^{-1} \bigr]^1_0,
\end{eqnarray}
where $\mathrm{Id}$ denotes the identity function on $\mathbb{R}$.
\end{proposition}

\begin{pf}
$\mathfrak{f}^\varepsilon_{T_{i\to o}}+\varepsilon(2\mathrm{Id}-1)$ increases
continuously from $\mathbb{R}$ onto $\mathbb{R}$, so $\{\mathfrak
{f}^\varepsilon_{T_{i\to
o}}+\varepsilon(2\mathrm{Id}-1)\}^{-1}$ exists and increases continuously
from $\mathbb{R}$ onto $\mathbb{R}$. Consequently,
the function $g\dvtx\mathbb{R}\to\mathbb{R}$ appearing in the
right-hand side of (\ref
{eqreceps}) is continuously increasing from $\mathbb{R}$ onto $\mathbb
{R}$, hence
invertible. Given $x \in\mathbb{R}$, it now remains to prove that
$t:=\mathfrak{f}
^{\varepsilon}_{(T,o)}(x)$ satisfies $g(t)=x$. By definition,
\begin{eqnarray}
\label{eq00} t & = & x+\partial\theta(o),
\end{eqnarray}
where $\theta$ denotes the $\varepsilon$-balanced allocation on $T$
with baseload $x$ at $o$ and $0$ elsewhere. Now fix $i\sim o$. The
restriction of $\theta$ to ${T_{i\to o}}$ is clearly an $\varepsilon
$-balanced allocation on ${T_{i\to o}}$ with baseload $\theta(o,i)$ at
$i$ and $0$ elsewhere. This is precisely the allocation appearing in
the definition of $\mathfrak{f}^\varepsilon_{T_{i\to o}}(\theta
(o,i))$, hence
\begin{eqnarray*}
\mathfrak{f}^\varepsilon_{T_{i\to o}}\bigl(\theta(o,i)\bigr) & = &
\partial \theta(i).
\end{eqnarray*}
Thus, the fact that $\theta$ is $\varepsilon$-balanced along $(o,i)$
may now be rewritten as
\begin{eqnarray}
\label{eq01} \theta(o,i) & = & \biggl[\frac{1}{2}+\frac{t-\mathfrak
{f}^\varepsilon_{T_{i\to
o}}(\theta(o,i))}{2\varepsilon}
\biggr]^1_0.
\end{eqnarray}
But by definition, $x_i:=\{\mathfrak{f}^\varepsilon_{T_{i\to
o}}+\varepsilon
(2\mathrm{Id}-1)\}^{-1}(t)$ is the unique solution to
\begin{eqnarray}
\label{eq02} x_i & = & \frac{1}{2}+\frac{t-\mathfrak{f}^\varepsilon_{T_{i\to
o}}(x_i)}{2\varepsilon}.
\end{eqnarray}
Comparing (\ref{eq01}) and (\ref{eq02}), we see that $\theta
(o,i)= [x_i ]^1_0$, that is, $\theta(i,o)=
[1-x_i
]^1_0$. Re-injecting this into (\ref{eq00}), we arrive exactly at the
desired $x=g(t)$.
\end{pf}

In the remainder of this section, we fix a vanishing sequence $\{
\varepsilon_n\}
_{n\geq1}$ and study the pointwise limit $\mathfrak{f}=\lim_{n\to
\infty}\mathfrak{f}^{\varepsilon
_n}_{(T,o)}$, when it exists. Note that $\mathfrak{f}$ needs not be invertible.
However, (\ref{eqnondnone}) and (\ref{eqcoercive}) guarantee that
$\mathfrak{f}
$ is nondecreasing with $\mathfrak{f}(\pm\infty)=\pm\infty$, so
that is admits a
well-defined right-continuous inverse
\begin{eqnarray*}
\label{eqrcinv} \mathfrak{f}^{-1}(t)&:=&\sup\bigl\{x\in\mathbb{R}\dvtx
\mathfrak {f}(x)\leq t\bigr\}, \qquad t\in\mathbb{R}.
\end{eqnarray*}

\begin{proposition}
\label{prcavity}
Assume that $\ell_o:=\lim_{n\to\infty}\partial\Theta_{\varepsilon
_n}(T,o)$
exists for each $o\in V$. Then
$\mathfrak{f}_{T_{i\to j}}:=\lim_{n\to\infty}\mathfrak
{f}^{\varepsilon_n}_{T_{i\to j}}$
exists pointwise for each $(i,j)\in\vec{E}$, and
\begin{eqnarray}
\label{eqreczero} \mathfrak{f}^{-1}_{T_{i\to j}}(t) & = & t-\sum
_{k\sim i, k\neq
j} \bigl[1-\mathfrak{f} ^{-1}_{T_{k\to i}}(t)
\bigr]^1_0,
\end{eqnarray}
for every $t\in\mathbb{R}$. Moreover, for every $o\in V$,
\begin{eqnarray}
\label{eqroot} \ell_o>t \quad & \Longleftrightarrow & \quad \sum
_{i\sim o} \bigl[1-\mathfrak {f}^{-1}_{T_{i\to
o}}(t)
\bigr]^1_0>t.
\end{eqnarray}
\end{proposition}

\begin{pf}
Fix $(i,j)\in\vec{E},x\in\mathbb{R}$, and let us show that $\{
\mathfrak{f}^{\varepsilon
_n}_{T_{i\to j}}(x)\}_{n\geq1}$ converges. By definition, $\mathfrak{f}
^{\varepsilon}_{{T_{i\to j}}}(x)=x+\partial\theta_\varepsilon(i)$,
where $\theta
_\varepsilon$ is the $\varepsilon$-balanced allocation on $T_{i\to j}$
with baseload $x$ at $i$ and $0$ elsewhere. Since the set of
allocations on $T_{i\to j}$ is compact, it is enough to consider two
subsequential limits $\theta,\theta'$ of $\{\theta_{\varepsilon_n}\}
_{n\geq1}$ and prove that $\partial\theta=\partial\theta'$.
Passing to
the limit in (\ref{eqbalb}), we know that $\theta,\theta'$ are
balanced with respect to the above baseload. Writing $V_{i\to j}$ for
the vertex set of $T_{i\to j}$, Proposition~\ref{lmdicho} reduces our
task to proving
\begin{eqnarray}
\label{eqtoshowl1}
\bigl\|\partial\theta-\partial\theta'\bigr\|_{\ell^1(V_{i\to j})} &
< & \infty.
\end{eqnarray}
Let ${\theta}^\star_\varepsilon$ be the restriction of $\Theta
_\varepsilon$ to $T_{i\to
j}$. Thus, $\theta^\star_\varepsilon$ is an allocation on $T_{i\to
j}$ and it is
$\varepsilon$-balanced with baseload $\theta^\star_\varepsilon
(j,i)$ at $i$ and $0$
elsewhere. Consequently, Proposition~\ref{prnonexp} guarantees that
for any finite $K\subseteq V_{i\to j}\setminus\{i\}$,
\begin{eqnarray*}
\bigl\|\partial{\theta}_\varepsilon-\partial{\theta}^\star_{\varepsilon}\bigr\|
_{\ell
^1(K)} & \leq& |x|+1.
\end{eqnarray*}
Applying this to $\varepsilon,\varepsilon'>0$ and using the triangle
inequality, we obtain
\[
\| \partial{\theta}_{\varepsilon}- \partial{\theta}_{\varepsilon'}\|_{\ell^1(K)}
 \leq  2|x|+2+\bigl\| \partial {\theta}^{\star}_{\varepsilon}-\partial {\theta}^{\star}_{\varepsilon'}\bigr\|_{\ell^1(K)}.
\]
Since $\{\partial\theta^{\star}_{\varepsilon_n}\}_{n\geq1}$ converges
by assumption, we may pass to the limit to obtain
$\|\partial\theta-\partial\theta'\|_{\ell^1(K)} \leq2|x|+2$. But $K$
is arbitrary, so (\ref{eqtoshowl1}) follows. This shows that
$\mathfrak{f}
_{T_{i\to j}}:=\lim_{n\to\infty}\mathfrak{f}^{\varepsilon
_n}_{T_{i\to j}}$ exists
pointwise. We now recall two classical facts about nondecreasing
functions $\mathfrak{f}\dvtx\mathbb{R}\to\mathbb{R}$ with
$\mathfrak{f}(\pm\infty)=\pm\infty$. First, $\mathfrak{f}
^{-1}$ is nondecreasing, so that its discontinuity set $\mathcal{D}(f^{-1})$
is countable. Second, the pointwise convergence $\mathfrak{f}_n\to
\mathfrak{f}$ implies $\mathfrak{f}
^{-1}_n(t)\to\mathfrak{f}^{-1}(t)$ for every $t\in\mathbb
{R}\setminus\mathcal{D}(\mathfrak{f}^{-1})$.
Consequently, letting $\varepsilon\to0$ in (\ref{eqreceps}) proves
(\ref{eqreczero}) for $t\notin\mathcal{D}:=\mathcal{D}(\mathfrak
{f}^{-1}_{T_{i\to j}})\cup\bigcup_{k\sim i} \mathcal{D}(\mathfrak{f}^{-1}_{T_{k\to i}})$. The equality
then extends to $\mathbb{R}
$ since $\mathcal{D}$ is countable and both sides of (\ref
{eqreczero}) are
right-continuous in $t$. Replacing $T_{i\to j}$ with $(T,o)$ in the
above argument shows that $\mathfrak{f}_{(T,o) }:=\lim_{n\to\infty
}\mathfrak{f}^{\varepsilon
_n}_{(T,o)}$ exists and satisfies
\begin{eqnarray*}
\mathfrak{f}^{-1}_{(T,o)}(t) & = & t-\sum
_{i\sim o} \bigl[1-\mathfrak {f}^{-1}_{T_{i\to
o}}(t)
\bigr]^1_0,  \qquad t\in\mathbb{R}.
\end{eqnarray*}
Finally, recall that $\mathfrak{f}^{\varepsilon
_n}_{(T,o)}(0)=\partial\Theta_{\varepsilon_n}(T,o)$
for all $n\geq1$, so that $\mathfrak{f}_{(T,o)}(0)=\ell_o$. But
$\mathfrak{f}
_{(T,o)}(0)>t\Longleftrightarrow\mathfrak{f}^{-1}_{(T,o)}(t)<0$ by
definition of
$\mathfrak{f}^{-1}_{(T,o)}$, so (\ref{eqroot}) follows.
\end{pf}

\section{Proof of Theorem~\texorpdfstring{\protect\ref{thmain2}}{2}}
\label{secthmain2}
In all this section, $t\in\mathbb{R}$ is fixed. We manipulate
networks rather
than graphs, where each $(i,j)\in\vec{E}$ is equipped with a mark
$\xi
(i,j)\in[0,1]$. The marks are assumed to satisfy the local recursion
\begin{equation}
\label{eqcavity}
\xi(i,j)= \biggl[1-t+\sum_{k\sim i,k\neq j}
\xi(k,i) \biggr]^1_0, \qquad (i,j)\in\vec{E}.
\end{equation}
We start with a simple lemma.

\begin{lemma}
\label{lmbreak}
$\partial\xi(i)\wedge\partial\xi(j)>t \Longleftrightarrow\xi
(i,j)+\xi(j,i)>1$.
\end{lemma}

\begin{pf}
We check the equivalence separately in each case. By assumption,
\begin{eqnarray}
\label{eqij} \xi(i,j) & = & \bigl[1-t+\partial\xi(i)-\xi(j,i)
\bigr]^1_0,
\\
\label{eqji} \xi(j,i) & = & \bigl[1-t+\partial\xi(j)-\xi(i,j)
\bigr]^1_0.
\end{eqnarray}
\begin{itemize}
\item{If ${0< \xi(i,j),\xi(j,i)<1}$}, then the equivalence trivially
holds since we may safely remove the truncation $[\cdot]^1_0$ from
(\ref
{eqij})--(\ref{eqji}) to obtain
\begin{eqnarray*}
\partial\xi(i)-t & = & \xi(i,j)+\xi(j,i)-1  =  \partial\xi(j)-t.
\end{eqnarray*}
\item{If ${\xi(j,i)=0}$}, then we have $1-t+\partial\xi(j)-\xi
(i,j)\leq0$ thanks to (\ref{eqji}), and hence $\partial\xi(j)\leq t$.
Thus, both sides of the equivalence are false.
\item{If ${\xi(i,j)=1,\xi(j,i)>0}$}, then using $\xi(i,j)=1$ in~(\ref
{eqij}) gives
$\partial\xi(i)-t\geq\xi(j,i)$ and since $\xi(j,i)>0$ we obtain
$\partial\xi(i)>t$. Similarly, using $\xi(j,i)>0$ in~(\ref{eqji})
gives $\partial\xi(j)>t+\xi(i,j)-1$ and since $\xi(i,j)=1$ we obtain
$\partial\xi(j)>t$. Thus, both sides of the equivalence are true.
\end{itemize}
The other possible cases follow by exchanging $\xi(i,j)$ and $\xi(j,i)$.
\end{pf}
We are ready for the proof of Theorem~\ref{thmain2}, which we divide
into two parts. The notation are those of Theorem~\ref{thmain2}, that
is, $\mu:=\operatorname{UGWT}(\pi)$, where $\pi$ is a fixed probability
distribution on $\mathbb{N}$ with finite, nonzero mean.

\begin{proposition}
If $Q\in\mathcal{P}([0,1])$ satisfies $Q=F_{\pi,t}(Q)$, then
\begin{eqnarray*}
\Psi_{\mathcal{L}_\mu}(t)
& \geq&
\frac{\mathbb{E}[D]}{2}\mathbb{P} (\xi_1+\xi _2>1 )-t\mathbb{P}
(\xi _1+\cdots+\xi_D>t ),
\end{eqnarray*}
where $D\sim\pi$ and where $\{\xi_k\}_{k\geq1}$ are i.i.d. with law
$Q$, independent of $D$.
\end{proposition}

\begin{pf}
Kolmogorov's extension theorem allows us to convert the consistency
equation $Q=F_{\pi,t}(Q)$ into a random rooted tree $\mathbb{T}\sim
\operatorname{UGWT}(\pi)$ equipped with marks satisfying (\ref{eqcavity}) a.s., such
that conditionally on the structure of $[\mathbb{T},o]_h$, the marks from
generation $h$ to $h-1$ are {i.i.d.} with law $Q$. This random
rooted network is easily checked to be unimodular. Thus, we may apply
Proposition~\ref{prdual} with $f=\mathbf{1}_{\partial\xi>t}$. By
Lemma~\ref
{lmbreak}, we have $\widehat{f}=\mathbf{1}_{\xi+\xi^*>1}$, and hence
\begin{eqnarray*}
\Psi_{\mathcal{L}_\mu}(t) & \geq& \tfrac{1}{2}\vec{\mu} \bigl(\xi +\xi^*>1
\bigr)-t\mu(\partial\xi>t).
\end{eqnarray*}
This is precisely the desired result, since we have by construction
\[
\mu(\partial\xi>t)=\mathbb{P} (\xi_1+\cdots+\xi_D>t ),\qquad
\vec{\mu } \bigl(\xi+\xi^*>1 \bigr)=\mathbb{E}[D]\mathbb{P} (\xi _1+
\xi_2>1 ),
\]
where $D\sim\pi$ and $\xi_1,\xi_2,\ldots$ are {i.i.d.} with law
$Q$, independent of $D$.
\end{pf}

\begin{proposition}
There exists $Q\in\mathcal{P}([0,1])$ with $Q=F_{\pi,t}(Q)$ and
\begin{eqnarray*}
\Psi_{\mathcal{L}_\mu}(t) &= & \frac{\mathbb{E}[D]}{2}\mathbb {P} (\xi_1+
\xi_2>1 )-t\mathbb{P} (\xi_1+\cdots+\xi_D >t
),
\end{eqnarray*}
where $D\sim\pi$ and where $\{\xi_k\}_{k\geq1}$ are {i.i.d.} with
law $Q$, independent of $D$.
\end{proposition}

\begin{pf}
Let ${\mathbb{T}}\sim\operatorname{UGWT}(\pi)$. Thanks to
Proposition~\ref{przero}, we have
\[
\partial\Theta_{\varepsilon}(\mathbb{T},o)\mathop{\longrightarrow
}_{\varepsilon\to0}^{L^2} \partial\Theta _0(\mathbb{T},o).
\]
In particular, there is a deterministic vanishing sequence $\varepsilon
_1,\varepsilon_2,\ldots$ along which the convergence holds
almost surely. This
almost-sure convergence automatically extends from the root to all
vertices, since under a unimodular measure $\mu$, \textit{everything
shows at the root} \cite{AldLyo07}, Lemma~2.3. More precisely,
\begin{eqnarray*}
\mu(A)=1 \quad& \Longrightarrow& \quad\mu(\widetilde{A})=1,
\end{eqnarray*}
for any Borel set $A\subseteq\mathcal{G}_{\star}$, where $\widetilde
{A}$ consists of
those $(G,o)\in\mathcal{G}_{\star}$ such that $(G,i)\in A$ for all
vertices $i$ of
$G$. Here, we apply it to $\mu=\operatorname{UGWT}(\pi)$ and
\[
A= \bigl\{(G,o)\in\mathcal{G}_{\star}\dvtx\partial\Theta
_{\varepsilon_n}(G,o) \mathop{\longrightarrow}_{n\to\infty}\partial
\Theta_{0}(G,o) \bigr\}.
\]
Thus, $\mathbb{T}$ satisfies almost surely the assumption of
Proposition~\ref
{prcavity}. Consequently, the marks $\xi(i,j) := [1-\mathfrak
{f}^{-1}_{\mathbb{T}
_{i\to j}}(t)]^1_0$ satisfy (\ref{eqcavity}) almost surely, and
\begin{eqnarray*}
\partial\Theta_0(\mathbb{T},o)>t & \quad\Longleftrightarrow\quad& \partial
\xi(o)>t.
\end{eqnarray*}
This ensures that $f=\mathbf{1}_{\partial\xi>t}$ satisfies the requirements
for equality in Proposition~\ref{prdual}, and we may then use
Lemma~\ref{lmbreak} to rewrite the conclusion as
\begin{eqnarray*}
\Psi_{\mathcal{L}_\mu}(t) & = & \tfrac{1}{2}\vec{\mu} \bigl(\xi +\xi^*>1
\bigr)-t\mu(\partial\xi>t).
\end{eqnarray*}
Now, $D=\operatorname{deg} (\mathbb{T},o)$ has law $\pi$ and conditionally on $D$, the
subtrees $\{\mathbb{T}_{i\to o}\}_{i\sim o}$ are {i.i.d.} copies
of a
homogenous Galton--Watson tree $\widehat{\mathbb{T}}$ with offspring
distribution $\widehat{\pi}$. Since $\xi(i,o)$ depends only on the
subtree $\mathbb{T}_{i\to o}$, we obtain
\[
\mu(\partial\xi>t)=\mathbb{P} (\xi_1+\cdots+\xi_D>t ),\qquad
\vec{\mu } \bigl(\xi+\xi^*>1 \bigr)=\mathbb{E}[D]\mathbb{P} (\xi _1+
\xi_2>1 ),
\]
where $\xi_1,\xi_2,\ldots$ are {i.i.d.} copies of $[1-\mathfrak{f}
^{-1}_{\widehat{\mathbb{T}}}(t)]^1_0$, independent of $D$. In turn, removing
the root of ${\widehat{\mathbb{T}}}$ splits it into a ${\widehat{\pi
}}$-distributed number of {i.i.d.} copies of $\widehat{\mathbb{T}}$, so
that the law $Q$ of $[1-\mathfrak{f}^{-1}_{\widehat{\mathbb{T}}}(t)]^1_0$ satisfies
$Q=F_{\pi,t}(Q)$.
\end{pf}

\section{Proof of Theorem~\texorpdfstring{\protect\ref{thmain3}}{3}}
\label{secthmain3}
In this final section, we prove Theorem~\ref{thmain3}. This main
ingredient is Proposition~\ref{prpairing}, which states that dense
subgraphs must be large under the pairing model. Fix a degree sequence
${\mathbf d}=\{d(i)\}_{1\leq i\leq n}$ and set $2m=\sum_{i=1}^nd(i)$. We
need two preparatory lemmas.

\begin{lemma}
\label{lmbin}
Fix a subset of vertices $S\subseteq\{1,\ldots,n\}$. Then the number of
edges of $\mathbb{G}[{\mathbf d}]$ with both end-points in $S$ is
stochastically
dominated by a binomial random variable with mean
$\frac{1}{m} (\sum_{i\in S} d_i )^2$.
\end{lemma}

\begin{pf}
We assume that $s:=\sum_{i\in S}d_i<m$, otherwise the claim is trivial.
It is classical that $\mathbb{G}[{\mathbf d}]$ can be generated
sequentially: at
each step $1\leq t\leq m$, a half-edge is selected and paired with a
uniformly chosen other half-edge. The selection rule is arbitrary, and
we choose to give priority to half-edges whose end-point lies in $S$.
Let $X_t$ be the number of edges with both end-points in $S$ after $t$
steps. Then $\{X_t\}_{0\leq t\leq m}$ is a Markov chain with $X_0=0$
and transitions
\[
X_{t+1} := %
\cases{ \displaystyle X_{t}+1, & $
\displaystyle\quad\mbox{with conditional probability }\frac{(s-X_{t}-t-1)^+}{2m-2t-1}$,\vspace
*{2pt}
\cr
X_{t}, & \quad$\mbox{otherwise}$.}
\]
For every $0\leq t < m$, the fact that $X_t\geq0$ ensures that
\begin{eqnarray*}
\frac{ (s-X_{t}-t-1 )^+}{2m-2t-1} & \leq& \frac
{s-t-1}{2m-2t-1}{\mathbf{1}}_{(t<s)}  \leq
 \frac{s}{2m}{\mathbf{1}}_{(t<s)},
\end{eqnarray*}
where the second inequality uses the condition $s< m$. This shows that
$X_m$ is in fact stochastically dominated by a binomial $(s,\frac{s}{2m})$, which is enough.
\end{pf}

\begin{lemma}
\label{lmcount}
Let $X_{k,r}$ be the number of induced subgraphs with $k$ vertices and
at least $r$ edges in $\mathbb{G}[{\mathbf d}]$. Then, for any $\theta>0$,
\begin{eqnarray*}
\mathbb{E} [X_{k,r} ] & \leq& \biggl(\frac{2r}{\theta^2m}
\biggr)^r \Biggl(\frac
{e}{k}\sum_{i=1}^ne^{\theta d_i}
\Biggr)^k. 
\end{eqnarray*}
\end{lemma}

\begin{pf}
First observe that if $Z\sim\operatorname{Bin}(n,p)$ then by a simple union-bound,
\begin{eqnarray*}
\mathbb{P} (Z\geq r ) & \leq& \pmatrix{n \cr r} p^r  \leq
\frac
{n^rp^r}{r!} =  \frac{\mathbb{E}[Z]^r}{r!}.
\end{eqnarray*}
Thus, Lemma~\ref{lmbin} ensures that the number $Z_S$ of edges with
both end-points in $S$ satisfies
\begin{eqnarray*}
\mathbb{P} (Z_S\geq r ) & \leq& \frac{1}{r! m^r} \biggl(\sum
_{i\in
S}d_i \biggr)^{2r}
\leq \biggl(\frac{2r}{\theta^{2}m} \biggr)^r\prod
_{i\in
S}e^{\theta d_i},
\end{eqnarray*}
where we have used the crude bounds $x^{2r}\leq(2r)!e^x$ and
$(2r)!/r!\leq(2r)^r$.
The result follows by summing over all $S$ with $|S|=k$ and observing that
\[
\sum_{|S|=k}\prod_{i\in S}e^{\theta d_i}
\leq\frac{1}{k!} \Biggl(\sum_{i=1}^ne^{\theta d_i}
\Biggr)^k\leq \Biggl(\frac{k}{e}\sum
_{i=1}^ne^{\theta d_i} \Biggr)^k.
\]
The second inequality follows from the classical lower-bound $k!\geq
 (\frac{k}{e} )^k$.
\end{pf}
We now fix $\{{\mathbf d_n}\}_{n\geq1}$ as in Theorem~\ref{thmain3}.
Let $Z^{(n)}_{\delta,t}$ be the number of subsets $\varnothing
\subsetneq
S\subseteq\{1,\ldots,n\}$ such that $|S|\leq\delta n$ and
$|E(S)|\geq
t |S|$ in $\mathbb{G}_n:=\mathbb{G}[{\mathbf d_n}]$.

\begin{proposition}
\label{prpairing}
For each $t>1$, there is $\delta>0$ and $\kappa<\infty$ such that
\[
\mathbb{E} \bigl[Z^{(n)}_{\delta,t} \bigr]\leq\kappa \biggl(
\frac
{\ln n}{n} \biggr)^{t-1},
\]
uniformly in $n\geq1$. In particular, $Z^{(n)}_{\delta,t}= 0$ w.h.p.
as $n\to\infty$.
\end{proposition}

\begin{pf}
The assumptions of Theorem~\ref{thmain3} guarantee that
for some $\theta>0$,
\[
\alpha:=\inf_{n\geq1} \Biggl\{\frac{1}{n}\sum
_{i=1}^n d_n(i) \Biggr\}>0 \quad
\mbox{and} \quad \lambda:=\sup_{n\geq1} \biggl\{\frac{1}{n}\sum
_{i\in
V}e^{\theta d_n(i)} \biggr\}<\infty.
\]
Now, fix $t>1$ and choose $\delta>0$ small enough so that $f(\delta
)<1$, where
\[
f(\delta)  :=  \biggl(1\vee\frac{2(1+t)}{\alpha\theta^2}
\biggr)^{t+1}e\lambda\delta^{t-1}.
\]
Using Lemma~\ref{lmcount} and the trivial inequality $kt\leq\lceil
kt\rceil\leq k(t+1)$, we have
\[
\mathbb{E} \bigl[X^{(n)}_{k,\lceil kt\rceil} \bigr]  \leq \biggl(
\frac
{2\lceil kt\rceil}{\theta^2k\alpha} \biggr)^{\lceil kt\rceil
}(e\lambda )^k \biggl(
\frac{k}{n} \biggr)^{\lceil kt\rceil-k} \leq f^k \biggl(
\frac{k} n \biggr).
\]
Since $f$ is increasing, we see that for any $1\leq m \leq\delta n$,
\begin{eqnarray*}
\mathbb{E} \bigl[Z^{(n)}_{\delta,t} \bigr]  =  \sum
_{k=1}^{\lfloor\delta
n\rfloor}\mathbb{E} \bigl[X^{(n)}_{k,\lceil kt\rceil}
\bigr] & \leq& \sum_{k=1}^{m-1}f^k
\biggl(\frac{m}{n} \biggr) + \sum_{k=m}^{\lfloor\delta n\rfloor}f^k({
\delta})
\\
& \leq& \frac{f({m}/{n})}{1-f({m}/{n})} + \frac{f(\delta
)^{m}}{1-f({\delta})}.
\end{eqnarray*}
Choose $m\sim c\ln n$ with $c$ fixed. As $n\to\infty$, the first term
is of order $(\frac{\ln n}{n})^{t-1}$ while the second is of order
$f(\delta)^{c\ln n} \ll (\frac{\ln n}{n})^{t-1}$, if $c$ is large enough.
\end{pf}

\begin{pf*}{Proof of Theorem~\protect\ref{thmain3}}
The assumptions on $\{{\mathbf d_n}\}_{n\geq1}$ are more than sufficient
to guarantee that a.s., the local weak limit of $\{\mathbb{G}_n\}
_{n\geq1}$ is
$\mu:=\operatorname{UGWT}(\pi)$ (see, e.g., \cite{Bor12}). Thus, the weak
convergence $\mathcal{L}_{\mathbb{G}_n}\to\mathcal{L}_\mu$ holds
a.s., by Theorem~\ref
{thmain1}. Now, if $t<\varrho(\mu)$ then $\mathcal{L}_\mu
((t,\infty)
)>0$, so the Portmanteau theorem ensures that $\liminf_{n}\mathcal
{L}_{\mathbb{G}
_n} ((t,\infty) )>0$ a.s.
Consequently,
\[
\mathbb{P} \bigl(\varrho(\mathbb{G}_n)\leq t \bigr)  =  \mathbb
{P} \bigl(\mathcal{L}_{\mathbb{G}_n} \bigl((t,\infty) \bigr)=0 \bigr) \mathop{
\longrightarrow}_{n\to\infty} 0.
\]
On the other-hand, if $t>\varrho(\mu)$ then $\mathcal{L}_\mu
([t,\infty
) )=0$, so the Portmanteau theorem gives $\mathcal{L}_{\mathbb
{G}_n}
((t,\infty) )\to0$ a.s. Thus, with $\delta$ as in
Proposition~\ref
{prpairing},
\[
\mathbb{P} \bigl(\varrho(\mathbb{G}_n)> t \bigr)  \leq  \mathbb {P}
\bigl(\mathcal{L}_{\mathbb{G}_n} \bigl([t,\infty) \bigr)>\delta \bigr)+
\mathbb{P} \bigl(Z^{(n)}_{\delta
,t}> 0 \bigr)\mathop{
\longrightarrow}_{n\to\infty}0.
\]
Note that the requirement $t>1$ is fulfilled, since $\varrho(\mu)\geq
1$. Indeed, every node in a tree of size $n$ has load $1-\frac{1}{n}$,
and the assumption $\pi_0+\pi_1<1$ guarantees that the size of the
random tree $\mathbb{T}\sim\operatorname{UGWT}(\pi)$ is unbounded.
\end{pf*}

\section*{Acknowledgements}
The work was initiated at the follow-up meeting of the
Newton Institute programme N$^{\mathrm{o}} 86$: \textit{Stochastic
Processes in
Communication Sciences}. The authors thank M. Lelarge, R. Sundaresan
and an anonymous referee for fruitful discussions.

%





\printaddresses

\begin{thebibliography}{28}

\bibitem{AgaMee01}
%
\begin{bbook}[mr]
\bauthor{\bsnm{Agarwal},~\bfnm{Ravi~P.}\binits{R.~P.}},
\bauthor{\bsnm{Meehan},~\bfnm{Maria}\binits{M.}} \AND
\bauthor{\bsnm{O'Regan},~\bfnm{Donal}\binits{D.}}
(\byear{2001}).
\btitle{Fixed Point Theory and Applications}.
\bseries{Cambridge Tracts in Mathematics}
\bvolume{141}.
\bpublisher{Cambridge Univ. Press},
\blocation{Cambridge}.
\bid{doi={10.1017/CBO9780511543005}, mr={1825411}}
\end{bbook}
%

\bptok{imsref}%
\endbibitem

\bibitem{AldLyo07}
%
\begin{barticle}[mr]
\bauthor{\bsnm{Aldous},~\bfnm{David}\binits{D.}} \AND
\bauthor{\bsnm{Lyons},~\bfnm{Russell}\binits{R.}}
(\byear{2007}).
\btitle{Processes on unimodular random networks}.
\bjournal{Electron. J.~Probab.}
\bvolume{12}
\bpages{1454--1508}.
\bid{doi={10.1214/EJP.v12-463}, issn={1083-6489}, mr={2354165}}
\end{barticle}
%

\bptok{imsref}%
\endbibitem

\bibitem{AldSte04}
%
\begin{bincollection}[mr]
\bauthor{\bsnm{Aldous},~\bfnm{David}\binits{D.}} \AND
\bauthor{\bsnm{Steele},~\bfnm{J.~Michael}\binits{J.~M.}}
(\byear{2004}).
\btitle{The objective method: Probabilistic combinatorial optimization
and local weak convergence}.
In \bbooktitle{Probability on Discrete Structures}.
\bseries{Encyclopaedia Math. Sci.}
\bvolume{110}
\bpages{1--72}.
\bpublisher{Springer},
\blocation{Berlin}.
\bid{doi={10.1007/978-3-662-09444-0_1}, mr={2023650}}
\end{bincollection}
%

\bptok{imsref}%
\endbibitem

\bibitem{Ald01}
%
\begin{barticle}[mr]
\bauthor{\bsnm{Aldous},~\bfnm{David~J.}\binits{D.~J.}}
(\byear{2001}).
\btitle{The {$\zeta(2)$} limit in the random assignment problem}.
\bjournal{Random Structures Algorithms}
\bvolume{18}
\bpages{381--418}.
\bid{doi={10.1002/rsa.1015}, issn={1042-9832}, mr={1839499}}
\end{barticle}
%

\bptok{imsref}%
\endbibitem

\bibitem{AldBan05}
%
\begin{barticle}[mr]
\bauthor{\bsnm{Aldous},~\bfnm{David~J.}\binits{D.~J.}} \AND
\bauthor{\bsnm{Bandyopadhyay},~\bfnm{Antar}\binits{A.}}
(\byear{2005}).
\btitle{A survey of max-type recursive distributional equations}.
\bjournal{Ann. Appl. Probab.}
\bvolume{15}
\bpages{1047--1110}.
\bid{doi={10.1214/105051605000000142}, issn={1050-5164}, mr={2134098}}
\end{barticle}
%

\bptok{imsref}%
\endbibitem

\bibitem{BenScr01}
%
\begin{barticle}[mr]
\bauthor{\bsnm{Benjamini},~\bfnm{Itai}\binits{I.}} \AND
\bauthor{\bsnm{Schramm},~\bfnm{Oded}\binits{O.}}
(\byear{2001}).
\btitle{Recurrence of distributional limits of finite planar graphs}.
\bjournal{Electron. J. Probab.}
\bvolume{6}
\bpages{no. 23, 13 pp. (electronic)}.
\bid{doi={10.1214/EJP.v6-96}, issn={1083-6489}, mr={1873300}}
\end{barticle}
%

\bptok{imsref}%
\endbibitem

\bibitem{Bil99}
%
\begin{bbook}[mr]
\bauthor{\bsnm{Billingsley},~\bfnm{Patrick}\binits{P.}}
(\byear{1999}).
\btitle{Convergence of Probability Measures},
\bedition{2nd} ed.
\bpublisher{Wiley},
\blocation{New York}.
\bid{doi={10.1002/9780470316962}, mr={1700749}}
\end{bbook}
%

\bptok{imsref}%
\endbibitem

\bibitem{Bol80}
%
\begin{barticle}[mr]
\bauthor{\bsnm{Bollob{\'a}s},~\bfnm{B{\'e}la}\binits{B.}}
(\byear{1980}).
\btitle{A probabilistic proof of an asymptotic formula for the number
of labelled regular graphs}.
\bjournal{European J. Combin.}
\bvolume{1}
\bpages{311--316}.
\bid{doi={10.1016/S0195-6698(80)80030-8}, issn={0195-6698}, mr={0595929}}
\end{barticle}
%

\bptok{imsref}%
\endbibitem

\bibitem{Bor12}
%
\begin{bmisc}[auto:parserefs-M02]
\bauthor{\bsnm{Bordenave},~\bfnm{Charles}\binits{C.}}
(\byear{2012}).
\bhowpublished{Lecture notes on random graphs and probabilistic
combinatorial optimization. Available at
\surl{http://www.math.univ-toulouse.fr/\textasciitilde bordenave/coursRG.pdf}.}
\end{bmisc}
%

\bptok{imsref}%
\endbibitem

\bibitem{BorLel13}
%
\begin{barticle}[mr]
\bauthor{\bsnm{Bordenave},~\bfnm{Charles}\binits{C.}},
\bauthor{\bsnm{Lelarge},~\bfnm{Marc}\binits{M.}} \AND
\bauthor{\bsnm{Salez},~\bfnm{Justin}\binits{J.}}
(\byear{2013}).
\btitle{Matchings on infinite graphs}.
\bjournal{Probab. Theory Related Fields}
\bvolume{157}
\bpages{183--208}.
\bid{doi={10.1007/s00440-012-0453-0}, issn={0178-8051}, mr={3101844}}
\end{barticle}
%

\bptok{imsref}%
\endbibitem

\bibitem{CaiSan07}
%
\begin{binproceedings}[mr]
\bauthor{\bsnm{Cain},~\bfnm{Julie~Anne}\binits{J.~A.}},
\bauthor{\bsnm{Sanders},~\bfnm{Peter}\binits{P.}} \AND
\bauthor{\bsnm{Wormald},~\bfnm{Nick}\binits{N.}}
(\byear{2007}).
\btitle{The random graph threshold for {$k$}-orientability and a fast
algorithm for optimal multiple-choice allocation}.
In \bbooktitle{Proceedings of the {E}ighteenth {A}nnual ACM-{SIAM}
{S}ymposium on {D}iscrete {A}lgorithms}
\bpages{469--476}.
\bpublisher{ACM},
\blocation{New York}.
\bid{mr={2482873}}
\end{binproceedings}
%

\bptok{imsref}%
\endbibitem

\bibitem{FerRam07}
%
\begin{binproceedings}[mr]
\bauthor{\bsnm{Fernholz},~\bfnm{Daniel}\binits{D.}} \AND
\bauthor{\bsnm{Ramachandran},~\bfnm{Vijaya}\binits{V.}}
(\byear{2007}).
\btitle{The {$k$}-orientability thresholds for {$G\sb{n,p}$}}.
In \bbooktitle{Proceedings of the {E}ighteenth {A}nnual ACM-{SIAM}
{S}ymposium on {D}iscrete {A}lgorithms}
\bpages{459--468}.
\bpublisher{ACM},
\blocation{New York}.
\bid{mr={2482872}}
\end{binproceedings}
%

\bptok{imsref}%
\endbibitem

\bibitem{FouKho13}
%
\begin{binproceedings}[mr]
\bauthor{\bsnm{Fountoulakis},~\bfnm{Nikolaos}\binits{N.}},
\bauthor{\bsnm{Khosla},~\bfnm{Megha}\binits{M.}} \AND
\bauthor{\bsnm{Panagiotou},~\bfnm{Konstantinos}\binits{K.}}
(\byear{2011}).
\btitle{The multiple-orientability thresholds for random hypergraphs}.
In \bbooktitle{Proceedings of the {T}wenty-{S}econd {A}nnual
ACM-{SIAM} {S}ymposium on {D}iscrete {A}lgorithms}
\bpages{1222--1236}.
\bpublisher{SIAM},
\blocation{Philadelphia, PA}.
\bid{mr={2858395}}
\bptnote{check year}%
\end{binproceedings}
%

\bptok{imsref}%
\endbibitem

\bibitem{GamNow06}
%
\begin{barticle}[mr]
\bauthor{\bsnm{Gamarnik},~\bfnm{David}\binits{D.}},
\bauthor{\bsnm{Nowicki},~\bfnm{Tomasz}\binits{T.}} \AND
\bauthor{\bsnm{Swirszcz},~\bfnm{Grzegorz}\binits{G.}}
(\byear{2006}).
\btitle{Maximum weight independent sets and matchings in sparse random
graphs. {E}xact results using the local weak convergence method}.
\bjournal{Random Structures Algorithms}
\bvolume{28}
\bpages{76--106}.
\bid{doi={10.1002/rsa.20072}, issn={1042-9832}, mr={2187483}}
\end{barticle}
%

\bptok{imsref}%
\endbibitem

\bibitem{GaoWor10}
%
\begin{bincollection}[mr]
\bauthor{\bsnm{Gao},~\bfnm{Pu}\binits{P.}} \AND
\bauthor{\bsnm{Wormald},~\bfnm{Nicholas~C.}\binits{N.~C.}}
(\byear{2010}).
\btitle{Load balancing and orientability thresholds for random
hypergraphs [extended abstract]}.
In \bbooktitle{S{TOC}'10---{P}roceedings of the 2010 ACM
{I}nternational {S}ymposium on {T}heory of {C}omputing}
\bpages{97--103}.
\bpublisher{ACM},
\blocation{New York}.
\bid{mr={2743258}}
\end{bincollection}
%

\bptok{imsref}%
\endbibitem

\bibitem{Haj90}
%
\begin{barticle}[mr]
\bauthor{\bsnm{Hajek},~\bfnm{Bruce}\binits{B.}}
(\byear{1990}).
\btitle{Performance of global load balancing by local adjustment}.
\bjournal{IEEE Trans. Inform. Theory}
\bvolume{36}
\bpages{1398--1414}.
\bid{doi={10.1109/18.59935}, issn={0018-9448}, mr={1080823}}
\end{barticle}
%

\bptok{imsref}%
\endbibitem

\bibitem{Haj96}
%
\begin{barticle}[mr]
\bauthor{\bsnm{Hajek},~\bfnm{Bruce}\binits{B.}}
(\byear{1996}).
\btitle{Balanced loads in infinite networks}.
\bjournal{Ann. Appl. Probab.}
\bvolume{6}
\bpages{48--75}.
\bid{doi={10.1214/aoap/1034968065}, issn={1050-5164}, mr={1389831}}
\end{barticle}
%

\bptok{imsref}%
\endbibitem

\bibitem{Jan09}
%
\begin{barticle}[mr]
\bauthor{\bsnm{Janson},~\bfnm{Svante}\binits{S.}}
(\byear{2009}).
\btitle{The probability that a random multigraph is simple}.
\bjournal{Combin. Probab. Comput.}
\bvolume{18}
\bpages{205--225}.
\bid{doi={10.1017/S0963548308009644}, issn={0963-5483}, mr={2497380}}
\end{barticle}
%

\bptok{imsref}%
\endbibitem

\bibitem{KhaSun12}
%
\begin{barticle}[mr]
\bauthor{\bsnm{Khandwawala},~\bfnm{Mustafa}\binits{M.}} \AND
\bauthor{\bsnm{Sundaresan},~\bfnm{Rajesh}\binits{R.}}
(\byear{2014}).
\btitle{Belief propagation for optimal edge cover in the random
complete graph}.
\bjournal{Ann. Appl. Probab.}
\bvolume{24}
\bpages{2414--2454}.
\bid{doi={10.1214/13-AAP981}, issn={1050-5164}, mr={3262507}}
\bptnote{check year}%
\end{barticle}
%

\bptok{imsref}%
\endbibitem

\bibitem{LecLel12}
%
\begin{binproceedings}[mr]
\bauthor{\bsnm{Leconte},~\bfnm{M.}\binits{M.}},
\bauthor{\bsnm{Lelarge},~\bfnm{M.}\binits{M.}} \AND
\bauthor{\bsnm{Massouli{\'e}},~\bfnm{L.}\binits{L.}}
(\byear{2012}).
\btitle{Convergence of multivariate belief propagation, with
applications to cuckoo hashing and load balancing}.
In \bbooktitle{Proceedings of the {T}wenty-{F}ourth {A}nnual
ACM-{SIAM} {S}ymposium on {D}iscrete {A}lgorithms}
\bpages{35--46}.
\bpublisher{SIAM},
\blocation{Philadelphia, PA}.
\bid{mr={3185378}}
\end{binproceedings}
%

\bptok{imsref}%
\endbibitem

\bibitem{Lel12}
%
\begin{binproceedings}[mr]
\bauthor{\bsnm{Lelarge},~\bfnm{M.}\binits{M.}}
(\byear{2012}).
\btitle{A new approach to the orientation of random hypergraphs}.
In \bbooktitle{Proceedings of the {T}wenty-{T}hird {A}nnual ACM-{SIAM}
{S}ymposium on {D}iscrete {A}lgorithms}
\bpages{251--264}.
\bpublisher{ACM},
\blocation{New York}.
\bid{mr={3205213}}
\end{binproceedings}
%

\bptok{imsref}%
\endbibitem

\bibitem{Luc92}
%
\begin{bincollection}[mr]
\bauthor{\bsnm{{\L}uczak},~\bfnm{Tomasz}\binits{T.}}
(\byear{1992}).
\btitle{Sparse random graphs with a given degree sequence}.
In \bbooktitle{Random Graphs, {V}ol.~2 ({P}ozna\'n, 1989)}.
\bseries{Wiley-Intersci. Publ.}
\bpages{165--182}.
\bpublisher{Wiley},
\blocation{New York}.
\bid{mr={1166614}}
\end{bincollection}
%

\bptok{imsref}%
\endbibitem

\bibitem{Lyo05}
%
\begin{barticle}[mr]
\bauthor{\bsnm{Lyons},~\bfnm{Russell}\binits{R.}}
(\byear{2005}).
\btitle{Asymptotic enumeration of spanning trees}.
\bjournal{Combin. Probab. Comput.}
\bvolume{14}
\bpages{491--522}.
\bid{doi={10.1017/S096354830500684X}, issn={0963-5483}, mr={2160416}}
\end{barticle}
%

\bptok{imsref}%
\endbibitem

\bibitem{Sal13}
%
\begin{barticle}[mr]
\bauthor{\bsnm{Salez},~\bfnm{Justin}\binits{J.}}
(\byear{2013}).
\btitle{Weighted enumeration of spanning subgraphs in locally
tree-like graphs}.
\bjournal{Random Structures Algorithms}
\bvolume{43}
\bpages{377--397}.
\bid{doi={10.1002/rsa.20436}, issn={1042-9832}, mr={3094425}}
\end{barticle}
%

\bptok{imsref}%
\endbibitem

\bibitem{SalSha09}
%
\begin{barticle}[mr]
\bauthor{\bsnm{Salez},~\bfnm{Justin}\binits{J.}} \AND
\bauthor{\bsnm{Shah},~\bfnm{Devavrat}\binits{D.}}
(\byear{2009}).
\btitle{Belief propagation: An asymptotically optimal algorithm for
the random assignment problem}.
\bjournal{Math. Oper. Res.}
\bvolume{34}
\bpages{468--480}.
\bid{doi={10.1287/moor.1090.0380}, issn={0364-765X}, mr={2554069}}
\end{barticle}
%

\bptok{imsref}%
\endbibitem

\bibitem{ShaSha94}
%
\begin{bbook}[mr]
\bauthor{\bsnm{Shaked},~\bfnm{Moshe}\binits{M.}} \AND
\bauthor{\bsnm{Shanthikumar},~\bfnm{J.~George}\binits{J.~G.}}
(\byear{1994}).
\btitle{Stochastic Orders and Their Applications}.
\bpublisher{Academic Press, Inc.},
\blocation{Boston, MA}.
\bid{mr={1278322}}
\end{bbook}
%

\bptok{imsref}%
\endbibitem

\bibitem{Ste02}
%
\begin{bincollection}[mr]
\bauthor{\bsnm{Steele},~\bfnm{J.~Michael}\binits{J.~M.}}
(\byear{2002}).
\btitle{Minimal spanning trees for graphs with random edge lengths}.
In \bbooktitle{Mathematics and Computer Science, II ({V}ersailles, 2002)}.
\bseries{Trends Math.}
\bpages{223--245}.
\bpublisher{Birkh\"auser},
\blocation{Basel}.
\bid{mr={1940139}}
\end{bincollection}
%

\bptok{imsref}%
\endbibitem

\bibitem{Hof13}
%
\begin{bmisc}[auto:parserefs-M02]
\bauthor{\bsnm{van~der Hofstad},~\bfnm{Remco}\binits{R.}}
(\byear{2013}).
\bhowpublished{Random graphs and complex networks.
Available at \surl{http://www.win.tue.nl/\textasciitilde rhofstad/NotesRGCN.html}.}
\end{bmisc}
%

\bptok{imsref}%
\endbibitem
\end{thebibliography}
\end{document}